\documentclass{article}

\usepackage[english]{babel}

\usepackage[a4paper,top=3cm,bottom=3cm,left=3cm,right=3cm]{geometry}
\usepackage{amsmath}
\usepackage{amssymb}
\usepackage{amsfonts}
\usepackage{amsthm}
\usepackage{graphicx}
\usepackage{xcolor}
\usepackage[unicode=true,pdfusetitle,
 bookmarks=true,bookmarksnumbered=false,bookmarksopen=false,
 breaklinks=false,pdfborder={0 0 0},pdfborderstyle={},backref=false,colorlinks=true]{hyperref}
\usepackage{bm}

\newcommand{\dd}{\text{d}}
\newcommand{\R}{\mathbb{R}}
\newcommand{\rx}{\mathbf{Q}}
\newcommand{\ry}{\bar{\mathbf{q}}}
\newcommand{\ux}{\mathbf{u}}
\newcommand{\uy}[1]{\tilde{\mathsf{#1}}}

\newcommand{\s}{\mathbb{S}}
\newcommand{\xx}{\mathbf{x}}
\newcommand{\y}{\mathbf{y}}
\newcommand{\ys}{\mathsf{y}}
\newcommand{\xs}{\mathsf{x}}
\newcommand{\z}{\mathbf{z}}
\newcommand{\n}{\mathbf{n}}

\newcommand{\ka}{\kappa}

\newcommand{\Oo}{\mathcal{O}}
\newcommand{\zero}{\mathbf{0}}

\newcommand{\nn}{\mathbf{\bar n}}

\newcommand{\norm}[1]{\left\lvert #1 \right\rvert}

\theoremstyle{definition}
\newtheorem*{example*}{\protect\examplename}
\theoremstyle{plain}
\newtheorem{theorem}{\protect\theoremname}[section] 
\theoremstyle{plain}

\theoremstyle{plain}
\newtheorem{remark}[theorem]{\protect\remarkname}
\theoremstyle{plain}

\theoremstyle{plain}
\newtheorem{notation}[theorem]{\protect\notationname}

\numberwithin{equation}{section}

\makeatother

\providecommand{\examplename}{Example}
\providecommand{\lemmaname}{Lemma}
\providecommand{\theoremname}{Theorem}
\providecommand{\remarkname}{Remark}
\providecommand{\notationname}{Notation}
\providecommand{\propositionname}{Proposition}

\usepackage[
backend=biber,
style=numeric,
sorting=nyt,
giveninits=true, maxbibnames=99
]{biblatex}
\date{}

\begin{document}
\title{Corrected Trapezoidal Rule-IBIM for linearized Poisson-Boltzmann equation}
\author{Federico Izzo\footnote{Corresponding author}\ $^{,}$\footnote{Department of Mathematics, KTH Royal Institute of Technology, Stockholm, Sweden (\href{mailto:izzo@kth.se}{izzo@kth.se})}
\and Yimin Zhong\footnote{ Department of Mathematics and Statistics, Auburn University, Auburn AL, USA (\href{mailto:yimin.zhong@auburn.edu}{yimin.zhong@auburn.edu})}
\and Olof Runborg\footnote{Department of Mathematics, KTH Royal Institute of Technology, Stockholm, Sweden (\href{mailto:olofr@kth.se}{olofr@kth.se})}
\and Richard Tsai\footnote{ Department of Mathematics and Oden Institute for Computational Engineering and Sciences, The University of Texas at Austin, Austin TX, USA (\href{mailto:ytsai@math.utexas.edu}{ytsai@math.utexas.edu})}
}

\maketitle
\begin{abstract}
    In this paper, we solve the linearized Poisson-Boltzmann equation, used to model the electric potential of macromolecules in a solvent. We derive a corrected trapezoidal rule with improved accuracy for a boundary integral formulation of the linearized Poisson-Boltzmann equation. 
    More specifically, in contrast to the typical boundary integral formulations, the corrected trapezoidal rule is applied to integrate a system of compacted supported singular integrals using uniform Cartesian grids in $\mathbb{R}^3$, without explicit surface parameterization.  
    A Krylov method, accelerated by a fast multipole method, is used to invert the resulting linear system. We study the efficacy of the proposed method, and compare it to an existing, lower order method.
    We then apply the method to the computation of electrostatic potential of macromolecules immersed in solvent. 
    The solvent excluded surfaces, defined by a common approach, are merely piecewise smooth, and we study the effectiveness of the method for such surfaces.


\end{abstract}


\section{Introduction}
We propose to solve the linearized Poisson-Boltzmann equation using an Implicit Boundary Integral Method (IBIM), discretized by a corrected trapezoidal rule (CTR). 
The linearized Poisson-Boltzmann equation can be used to model the electrostatics of charged macromolecule-solvent systems.
Numerical modeling and simulation of such systems
is an active and relevant research topic. Electrostatic interactions with the solvent significantly affect the overall macromolecule behavior. They are relevant for example in electrochemistry, as an
aqueous solvent plays a significant role in the dynamical processes of biological molecules (see \cite{baker2005biomolecular,fogolari2002poisson} for comprehensive introductions to the problem). 

When describing the electrostatic interactions, implicit solvent methods approximate the problem by simplifying the description of the solvent environment and consequently greatly reducing the degrees of freedom. The frameworks for approximating implicitly the electrostatic interactions are for example the Coulomb-field approximation, the generalized Born models, and the Poisson-Boltzmann theory. The latter two are more popular, and among them the Poisson-Boltzmann theory offers a more detailed representation at the cost of more expensive computations (see also \cite{feig2004recent}). 

The Poisson-Boltzmann equation, which describes the electrostatic field determined by the interactions in the Poisson-Boltzmann theory, is nonlinear. The solution to its linearized form is however widely used, both because it represents a valid approximation to the nonlinear solution in certain settings, and because it can be used iteratively to find the original nonlinear solution. Hence, numerically solving the linearized Poisson-Boltzmann equation is still an active and relevant problem (see e.g. \cite{juffer1991electric,helgadottir2011poisson,lu2008recent,zhou2014variational,zhong2018implicit}).

In this work, we propose a fast numerical method for solving the linearized Poisson-Boltzmann equation on surfaces:
\begin{equation} \label{eq:linearized-Poisson-Boltzmann}
\begin{array}{rll}
    -\nabla \cdot \left( \varepsilon_I \nabla \psi(\xx) \right) & = \sum_{k=1}^{N_c} q_k \delta(\xx-\z_k), & \text{in }\Omega, \\[0.2cm]
    -\nabla \cdot \left( \varepsilon_E \nabla \psi(\xx) \right) & = -\bar\kappa^2_T \psi(\xx)  & \text{in }\bar\Omega^c ,\\[0.2cm]
    \left.\psi(\xx)\right|_{\Gamma^-} & = \left.\psi(\xx)\right|_{\Gamma^+} & \text{on }\Gamma,\\
    \varepsilon_I \left.\dfrac{\partial \psi}{\partial \n}(\xx)\right|_{\Gamma^-} & = \varepsilon_E \left.\dfrac{\partial \psi}{\partial \n}(\xx)\right|_{\Gamma^+} & \text{on }\Gamma,\\
    |\xx|\psi(\xx)\to 0, & |\xx|^2\left| \nabla \psi(\xx) \right|\to 0, & \text{as }|\xx|\to \infty.
\end{array}
\end{equation}
In \eqref{eq:linearized-Poisson-Boltzmann}, $\psi$ represents the electrostatic potential, $N_c$ is the number of atoms composing the macromolecules, which have centers $\{\z_k\}_{k=1}^{N_c}$, radii $\{r_k\}_{k=1}^{N_c}$, and charge numbers $\{q_k\}_{k=1}^{N_c}$ respectively. 
$\Gamma$ represents the closed surface which separates the region occupied by the macromolecule (here denoted by $\Omega$) and the rest of the space $\bar\Omega^c:=\R^3\setminus \bar\Omega=\R^3\setminus(\Omega \cup \Gamma)$.
We will describe in Section~\ref{sec:num:smooth} how we test the method on simple smooth surfaces $\Gamma$, and then in Section~\ref{sec:num:application} how we apply it to actual solvent-molecule interfaces.
The operator $\frac{\partial \psi}{\partial \n}(\xx)$ represents the normal derivative of $\psi$ in $\xx\in\Gamma$ with normal pointing outward from $\Gamma$. 
The parameters $\varepsilon_I$ and $\varepsilon_E$ are the dielectric constants inside and outside $\Omega$ respectively, and $\bar\kappa_T$ is a screening parameter. 
The radiation conditions on the last row are needed to ensure uniqueness of the solutions.

We use the following boundary integral formulation \cite{juffer1991electric} to find the solution to the equations \eqref{eq:linearized-Poisson-Boltzmann}:
\begin{align}
    \dfrac{1}{2}\left(1+\dfrac{\varepsilon_E}{\varepsilon_I}\right) \psi(\xx) &+ \int_\Gamma {K_{11}(\xx,\y)} \psi(\y) \dd\y - \int_\Gamma {K_{12}(\xx,\y)} \dfrac{\partial\psi}{\partial \n}(\y) \dd\y = \sum_{k=1}^{N_c} \dfrac{q_k}{\varepsilon_I} G_0(\xx,\z_k), \label{eq:PB-BIE-1}\\
    \dfrac{1}{2}\left(1+\dfrac{\varepsilon_I}{\varepsilon_E}\right) \dfrac{\partial\psi}{\partial \n}(\xx) &+ \int_\Gamma {K_{21}(\xx,\y)} \psi(\y) \dd\y - \int_\Gamma {K_{22}(\xx,\y)} \dfrac{\partial\psi}{\partial \n}(\y) \dd\y = \sum_{k=1}^{N_c} \dfrac{q_k}{\varepsilon_I} \dfrac{\partial G_0}{\partial \n_x}(\xx,\z_k). \label{eq:PB-BIE-2}
\end{align}
where
\begin{equation} \label{eq:Kij-kernels}
    \begin{array}{rlrl}
        K_{11}(\xx,\y) :=& \dfrac{\partial G_0}{\partial \n_y}(\xx,\y)-\dfrac{\varepsilon_E}{\varepsilon_I}\dfrac{\partial G_\ka}{\partial \n_y}(\xx,\y), & 
        K_{12}(\xx,\y) :=& {G_0}(\xx,\y)-{G_\ka}(\xx,\y), \\[0.4cm]
        K_{21}(\xx,\y) :=& \dfrac{\partial^2 G_0}{\partial \n_x \partial \n_y}(\xx,\y)-\dfrac{\partial^2 G_\ka}{\partial \n_x \partial \n_y}(\xx,\y), & 
        K_{22}(\xx,\y) :=& \dfrac{\partial G_0}{\partial \n_x}(\xx,\y)-\dfrac{\varepsilon_I}{\varepsilon_E}\dfrac{\partial G_\ka}{\partial \n_x}(\xx,\y).
    \end{array}
\end{equation}
The fundamental solutions which define the kernels in \eqref{eq:Kij-kernels} are
\begin{equation}\label{eq:layer-kernels}
    G_0(\xx,\y) := \dfrac{1}{4\pi|\xx-\y|}, \quad G_\kappa(\xx,\y) := \dfrac{e^{-\kappa|\xx-\y|}}{4\pi|\xx-\y|}.
\end{equation}
The integrals in \eqref{eq:PB-BIE-1}-\eqref{eq:PB-BIE-2} are well defined because of the integrability of the kernels for $\xx=\y$:
\begin{align} \label{eq:kij-asymptotics}
    K_{11}(\xx,\y) \sim \Oo(|\xx-\y|^{-1}), \quad & K_{21}(\xx,\y) \sim \Oo(|\xx-\y|^{-1}),\quad K_{22}(\xx,\y) \sim \Oo(|\xx-\y|^{-1}),\\
    & K_{12}(\xx,\y) \sim \Oo(1). \nonumber
\end{align}
As in \cite{zhong2018implicit}, we write these integrals in the Implicit Boundary Integral Method (IBIM) framework (see \cite{kublik2016integration, kublik2013implicit}). \\

We point out that the kernel $K_{21}$ is integrable even though it is defined as the difference between two hypersingular kernels which are not Cauchy integrable. 
This is the power of the boundary integral formulation derived by Juffer \emph{et al.} \cite{juffer1991electric}. 
The kernel $K_{21}$ proves to be the bottleneck of the accuracy order attainable given a second order approximation of the surface around the target point $\xx$. 
If the approximation is improved to third order, it is possible to improve the order of accuracy for the kernels $K_{11}$, $K_{22}$, and $K_{12}$ but not for $K_{21}$. 
Only a fourth order approximation allows us to increase the order of accuracy for $K_{21}$. 
This will be explained in more detail in Section~\ref{subsec:expansions}.

\section{Numerical solution of the volumetric boundary integral formulation }
\label{sec:BIE-CTR}

In this section we will provide the necessary background information for the numerical solution of \eqref{eq:PB-BIE-1}-\eqref{eq:PB-BIE-2}. 

In Section~\ref{subsec:ImplicitParametrization} we present the volumetric formulation of the surface integrals via non-parametrized surface representation. 
In Section~\ref{subsec:quadrature} we present the discretization \eqref{eq:corr-trapz-rule-sys1}-\eqref{eq:corr-trapz-rule-sys2} of the original BIEs \eqref{eq:PB-BIE-1}-\eqref{eq:PB-BIE-2} using the framework shown in the previous subsection.
In Section~\ref{subsec:Kreg} we briefly present the kernel regularization (K-reg) developed in \cite{zhong2018implicit}.
In Section~\ref{subsec:CorrTR} we present the formally second order accurate quadrature rule (CTR2) based on the trapezoidal rule which we use to approximate the singular volume integrals. 
In Section~\ref{subsec:overview-algorithm} we present an overview of the algorithms used to apply the quadrature rule given the surface $\Gamma$.
In Section~\ref{sec:num:smooth} we present numerical tests on smooth surfaces to show the order of accuracy of CTR2 and compare it to K-reg.

    \subsection{Extensions of the singular boundary integral operators}
    
\label{subsec:ImplicitParametrization}

Let $\Omega \subset\mathbb{R}^3$ be a bounded open set with $C^2$ boundaries, and $\partial\Omega=:\Gamma$. 
We shall refer to $\Gamma$ as the surface. Let $f$ be a function defined on $\Gamma$. 
In this section we present an approach for extending a boundary integral
\begin{equation}\label{eq:surfint}
    \int_{\Gamma}f(\y)\text{d}\sigma_{\y},
\end{equation}
to a volumetric integral around the surface.
Instead of parameterizations, this approach relies on the Euclidean distance to the surface, and its derivatives. More precisely, we define the \textit{signed distance function}
\begin{equation}\label{eq:signed-distance}
    d_{\Gamma}(\mathbf{x}):=\begin{cases}
\quad\min_{\mathbf{y}\in\Gamma}\left\|\mathbf{x}-\mathbf{y}\right\|, & \text{ if } \mathbf{x}\in \Omega\\
-\min_{\mathbf{y}\in\Gamma}\left\|\mathbf{x}-\mathbf{y}\right\|, & \text{ if } \mathbf{x}\in \Omega^c 
\end{cases}
\end{equation}
and the \textit{closest point projection} 
\begin{equation}\label{eq:closest-point-projection}
    P_{\Gamma}(\mathbf{x}):=\text{argmin}_{\mathbf{y}\in\Gamma}\left\|\mathbf{x}-\mathbf{y}\right\|.
\end{equation}
If there is more than one global minimum, we pick one randomly from the set. 
Let $\mathcal{C}_\Gamma$ denote the set of points in $\mathbb{R}^3$ which are equidistant to at least two distinct points on $\Gamma$.
The \emph{reach} $\tau$ is defined as $\inf_{\mathbf{x}\in \Gamma, \mathbf{y}\in\mathcal{C}_\Gamma} \|\mathbf{x}-\mathbf{y}\|$. Clearly, $\tau$ is restricted by the local geometry (the curvatures) and the global structure of $\Gamma$ (the Euclidean and geodesic distances between any two points on $\Gamma$). 

In this paper, we assume that $\Gamma$ is $C^2$ and thus has a non-zero reach. 
Let $T_{\varepsilon}$ denote the set of points of distance at most $\varepsilon$ from $\Gamma$:
\begin{equation}
    T_{\varepsilon}=\{\mathbf{x}\in\mathbb{R}^{3}:\ |d_{\Gamma}(\mathbf{x})|\leq \varepsilon\}.
\end{equation}
Then $P_\Gamma$ is a diffeomorphism between $\Gamma$ and the level sets of $d_\Gamma$ in $T_\epsilon$. Furthermore,
\[
P_\Gamma(\mathbf{x}) = \mathbf{x}-d_\Gamma(\mathbf{x})\nabla d_\Gamma(\mathbf{x}),~~~~\mathbf{x}\in T_\varepsilon.
\]

We define the \emph{extension}
of the integrand $f$ by
\begin{equation}\label{eq:extension-restriction-f}
    \overline{f}(\mathbf{x}) := f(P_\Gamma \mathbf{x}),~~~\mathbf{x}\in\mathbb{R}^3.
\end{equation} 
As in \cite{kublik2013implicit, kublik2016integration}, we can then rewrite the
surface integral (\ref{eq:surfint}):
\begin{equation}\label{eq:surface-integral-IBIM}
    \int_{\Gamma}f(\y)\text{d}\sigma_{\y} = \int_{\mathbb{R}^3}\overline{f}(\xx)\delta_{\Gamma,\varepsilon}(\y)\text{d}\y = \int_{T_\varepsilon}\overline{f}(\xx)\delta_{\Gamma,\varepsilon}(\y)\text{d}\y,
\end{equation}
where 
\begin{align}
    &\delta_{\Gamma,\varepsilon}(\y):=J_{\Gamma}(\y)\delta_{\varepsilon}(d_{\Gamma}(\y)), \quad \y\in \R^3,\nonumber\\
    &\delta_\varepsilon(\eta):= \frac{1}{\varepsilon}\phi\left(\frac{\eta}{\varepsilon}\right),\quad \phi\in C^\infty(\mathbb{R})\text{ supported in }[-1,1], \text{ and }\int_\R \phi(x)\dd x =1, \label{eq:phi-averaging-kernel}
\end{align}
and $J_{\Gamma}(\y')$ is the Jacobian of the transformation from  $\Gamma$ to the level set $\Gamma_{\eta}:=\{\y\in\mathbb{R}^{n}\ :\ d_{\Gamma}(\y)=\eta\}.$ In $\mathbb{R}^3$,
$J_\Gamma(\y')$ is a quadratic polynomial in $d_\Gamma(\y^\prime)$:
\begin{equation}\label{eq:jacobian}
    J_\Gamma(\y'):=1 + 2d_\Gamma(\y') \mathcal{H}(\y') + d_\Gamma(\y')^2 \mathcal{G}(\y').
\end{equation}
where $\mathcal{H}(\y')$ and $\mathcal{G}(\y')$ are respectively the mean and Gaussian curvatures of $\Gamma_\eta$ at $\y'$.
The curvatures as well as the corresponding principal directions can be extracted from $P_\Gamma$; see \cite{kublik2016integration} for more detail.

Our primary focus is when $f(\y)$ is replaced by a function of the kind $K(\xx,\y)\zeta(\y)$, 
corresponding to the kernels and potentials found in \eqref{eq:PB-BIE-1}-\eqref{eq:PB-BIE-2}, \eqref{eq:Kij-kernels}), and the integral operators 
of form
 \begin{equation} \label{eq:general_S} 
  \mathcal{J}_i[\zeta_1,\zeta_2](\xx) := \int_\Gamma K_{i1}(\xx,\y)\zeta_1(\y)\dd\y - \int_\Gamma  K_{i2}(\xx,\y)\zeta_2(\y) \dd\y, 
    \quad \xx \in \R^3,~~i=1,2.
\end{equation}
We can now write the operators \eqref{eq:general_S} in volumetric form:
\begin{equation}\label{eq:volume_potential_SL}
        \overline{\mathcal{J}}_i[\rho_1,\rho_2](\mathbf{x}) := \int_{T_\varepsilon} \overline{K}_{i1}(\mathbf{x},\mathbf{y})\rho_1(\mathbf{y}) \delta_{\Gamma, \varepsilon}(\mathbf{y}) \text{d}\mathbf{y} - \int_{T_\varepsilon} \overline{K}_{i2}(\mathbf{x},\mathbf{y})\rho_2(\mathbf{y}) \delta_{\Gamma, \varepsilon}(\mathbf{y}) \text{d}\mathbf{y},~~~\mathbf{x}\in T_\epsilon,~~~i=1,2,
\end{equation}
where 
\begin{equation}\label{eq:restriction_of_K}
    \overline{K} (\mathbf{x},\mathbf{y}):=K(\mathbf{x},P_\Gamma \mathbf{y})\ ,\ \ \mathbf{x},\mathbf{y}\in\R^n\,.
\end{equation}
It is important to notice that if $K(\mathbf{x},\mathbf{y})$ is singular for $\mathbf{x}=\mathbf{y}$, then $\overline{K}(P_\Gamma\mathbf{x},\mathbf{y})$ is singular on the set
\begin{equation}\nonumber
    \{ (\mathbf{x},\mathbf{y})\in \mathbb{R}^n\times\mathbb{R}^n: P_\Gamma \mathbf{x} = P_\Gamma \mathbf{y}\}, 
\end{equation}
i.e. for a fixed $\mathbf{x}^*\in \Gamma$, $\overline{K}(\mathbf{x}^*,\mathbf{y})$ is singular along the normal line passing through $\mathbf{x}^*$,
while $K(\mathbf{x}^*,\mathbf{y})$ is singular in a point.

Finally, with 
\begin{equation} \label{eq:PB-lambda-g}
    \begin{array}{rlrl}
        \lambda_1 & := \dfrac{1}{2}\left(1+\dfrac{\varepsilon_E}{\varepsilon_I}\right), & \lambda_2 & := \dfrac{1}{2}\left(1+\dfrac{\varepsilon_I}{\varepsilon_E}\right),  \\[0.4cm]
        g_1(\xx) & :=  \sum_{k=1}^{N_c} \dfrac{q_k}{\varepsilon_I} G_0(\xx,\z_k), & g_2(\xx) & := \sum_{k=1}^{N_c} \dfrac{q_k}{\varepsilon_I} \dfrac{\partial G_0}{\partial \n_x}(\xx,\z_k),
    \end{array}
\end{equation}
the volumetric forms of equations \eqref{eq:PB-BIE-1}-\eqref{eq:PB-BIE-2} are derived:
\begin{equation} \label{eq:IBIE}
        \lambda_i\rho_i(\xx) + \overline{\mathcal{J}}_i[\rho_1,\rho_2](P_\Gamma\xx) = g_i(P_\Gamma\xx),
    \quad \xx\in T_\varepsilon,~~i=1,2.
\end{equation}
The solutions $\rho_1$, $\rho_2$ will coincide with the constant extensions along the normals of $\psi$ and $\psi_n$ respectively: $\rho_1(\xx)\equiv \psi(P_\Gamma(\xx))$,  $\rho_2(\xx)\equiv \psi_n(P_\Gamma(\xx))$; see, for example, 
the arguments given in \cite{izzo2022corrected}.

    \subsection{Quadrature rules on uniform Cartesian grids }
    \label{subsec:quadrature}

In this paper, we derive numerical quadratures for the singular integral operators $\overline{\mathcal{J}}_i[\rho_1,\rho_2](\xx)$ for $\xx\in\Gamma$ (equivalently, $\overline{\mathcal{J}}_i[\rho_1,\rho_2](P_\Gamma\xx)$ for $\xx \in T_\varepsilon$) for $i=1,2$. 
The quadrature rules will be constructed based on the trapezoidal rule for the grid nodes $T^h_\varepsilon := T_\varepsilon \cap h\mathbb{Z}^3$, which corresponds to the portion of the uniform Cartesian grid $h\mathbb{Z}^3$ within $T_\varepsilon$.
Since the integrand in \eqref{eq:volume_potential_SL} is singular for $\mathbf{x}\in\Gamma$, the trapezoidal rule should be corrected near 
$\mathbf{x}$ for faster convergence. 
{Correction will be defined by summing the judiciously derived weights over a set of grid nodes denoted by $N_h(\mathbf{x}).$ The sum will be denoted by $\mathcal{R}_h(\mathbf{x})$. }
Ultimately, the quadrature for $\overline{\mathcal{J}}_i[\rho_1,\rho_2](P_\Gamma\mathbf{x})$ 
will involve the regular Riemann sum of the integrand in $T^h_\varepsilon\setminus N_h(\mathbf{x})$, and  the correction $\mathcal{R}_h(\mathbf{x})$ in $N_h(\mathbf{x})$:
\begin{align*}
  \overline{\mathcal{J}}_1[\rho_1,\rho_2](\xx) \approx & \
   h^3\sum_{\y_m \in T_\varepsilon^h \setminus N_h(\xx)} \overline{K}_{11}(\xx, \y_{m})\rho_1(\y_{m})\delta_{\Gamma, \varepsilon}(\y_m) + h^2\sum_{\y_m\in N_h(\xx)} \omega_m^{(11)} \rho_1(\y_{m})\delta_{\Gamma, \varepsilon}(\y_m)\\ 
   & - h^3\sum_{\y_m \in T_\varepsilon^h \setminus N_h(\xx)} \overline{K}_{12}(\xx, \y_{m})\rho_2(\y_{m})\delta_{\Gamma, \varepsilon}(\y_m) - h^3\sum_{\y_m\in N_h(\xx)} \omega_m^{(12)} \rho_2(\y_{m})\delta_{\Gamma, \varepsilon}(\y_m),\nonumber
\end{align*}
and
\begin{align*}
  \overline{\mathcal{J}}_2[\rho_1,\rho_2](\xx) \approx & \
   h^3\sum_{\y_m \in T_\varepsilon^h \setminus N_h(\xx)} \overline{K}_{21}(\xx, \y_{m})\rho_1(\y_{m})\delta_{\Gamma, \varepsilon}(\y_m) + h^2\sum_{\y_m\in N_h(\xx)} \omega_m^{(21)} \rho_1(\y_{m})\delta_{\Gamma, \varepsilon}(\y_m)\\ 
   & - h^3\sum_{\y_m \in T_\varepsilon^h \setminus N_h(\xx)} \overline{K}_{22}(\xx, \y_{m})\rho_2(\y_{m})\delta_{\Gamma, \varepsilon}(\y_m) - h^2\sum_{\y_m\in N_h(\xx)} \omega_m^{(22)} \rho_2(\y_{m})\delta_{\Gamma, \varepsilon}(\y_m),\nonumber 
\end{align*}
where $\omega_m^{(ij)}$ is a weight which depends on the kernel $K_{ij}$, on the grid node $\y_m$, and on the principal curvatures and directions of the surface in $\xx\in\Gamma$.

We now have all the ingredients to write out the linear system we need to solve: $\forall \y_k\in T_\varepsilon^h$,
\begin{align}
    \lambda_1 \rho_1(\y_k) +& h^3\sum_{\y_m \in T_\varepsilon^h\setminus N_h(\y_k)} \left(\overline{K}_{11}(P_\Gamma\y_k, \y_m )  \rho_1(\y_m) - \overline{K}_{12}(P_\Gamma\y_k,\y_m ) \rho_2(\y_m) \right) \delta_{\Gamma,\varepsilon}(\y_m) \nonumber\\
    +& h^2\sum_{\y_m\in N_h(\y_k)} \left( \omega^{(11)}_{mk} \rho_1(\y_m) - h\,\omega^{(12)}_{mk} \rho_2(\y_m) \right) \delta_{\Gamma,\varepsilon}(\y_m) = g_1(P_\Gamma\y_k), \label{eq:corr-trapz-rule-sys1}\\
    \lambda_2 \, \rho_2(\y_k) +& h^3\sum_{\y_m \in T_\varepsilon^h\setminus N_h(\y_k)} \left(\overline{K}_{21}(P_\Gamma\y_k,\y_m)  \rho_1(\y_m) - \overline{K}_{22}(P_\Gamma\y_k,\y_m)\rho_2(\y_m) \right) \delta_{\Gamma,\varepsilon}(\y_m) \nonumber\\
    +& h^2\sum_{\y_m \in N_h(\y_k)} \left( \omega^{(21)}_{mk} \rho_1(\y_m) - \omega^{(22)}_{mk} \rho_2(\y_m) \right) \delta_{\Gamma,\varepsilon}(\y_m) = g_2(P_\Gamma\y_k),\label{eq:corr-trapz-rule-sys2}
\end{align}
where $\rho_1(\y_k)\approx \bar\psi(\y_k)$ and $\rho_2(\y_k) \approx \bar\psi_n(\y_k)$. 
Corresponding to the target node $\y_k$ the set $N_h(\y_k)$ contains all the correction nodes, and $\omega_{mk}^{(ij)} = \omega[s^{(ij)}; \alpha_{mk},\beta_{mk}]$, $i,j=1,2$, are the correction weights, dependent on the kernel they correct\footnote{Please note that the weight for the kernel $K_{12}$ has a different scaling compared to the other three. This is because in order to reach order of accuracy two for such kernel (see \eqref{eq:kij-asymptotics}) the correction weight is zero, so we apply the correction necessary to reach order three.} via the function $s^{(ij)}$, and on  $(\alpha_{mk},\beta_{mk})$, parameters dependent on $h$, $\y_m$, and $\y_k$. The details of how the corrections nodes are chosen are presented in Section~\ref{subsec:CorrTR}, while the definition and computation of the weights is presented in Section~\ref{sub:weights}.

The contribution of this paper are two quadrature rules. One is a high order, trapezoidal rule-based, quadrature rule for $\mathcal{J}$ via $\overline{\mathcal{J}}$, effective for smooth surfaces which are globally $C^2$. The second is a hybrid rule which combines the previous trapezoidal rule-based with the constant regularization from \cite{zhong2018implicit}.

    \subsection{Kernels regularization in IBIM - K-reg}
    \label{subsec:Kreg}
The principle of the regularization proposed in \cite{zhong2018implicit} is to substitute the kernel $\overline{K}$ with a regularized version $\tilde{K}$ around the singularity point, so that the integral over the domain would be as close as possible to the original one. Thus the technique is independent of the quadrature rule. Given a parameter $\tau>0$ and $\xx,\y\in T_\varepsilon$, the regularization of the kernel is
\begin{equation}\label{eq:regularized-kernel}
  \overline{K}_{\tau}(\xx,\y) := 
  \begin{cases}
    \overline{K}(\xx,\y), & \text{if } \norm{\xx-\y}_P\geq \tau, \\
    C_{\Gamma,\tau}, & \text{if } \norm{\xx-\y}_P < \tau, \\
  \end{cases}
\end{equation}
where $\norm{\xx-\y}_P$ is the distance between the projections of $\xx$ and $\y$ on the tangent plane at $P_\Gamma\xx$, and $C_{\Gamma,\tau}$ is a constant dependent on the surface $\Gamma$ and on the parameter $\tau$.  
The constant $C_{\Gamma,\tau}$ is constructed so that it behaves in the same way as $\overline{K}$ in a neighborhood of $\xx$ dependent on $\tau$. Specifically, given $V(P_\Gamma\xx;\tau)$ disc of radius $\tau$ on the tangent plane of $\Gamma$ at $P_\Gamma\xx$, it is defined as
\[
  C_{\Gamma,\tau} = \dfrac{1}{V(P_\Gamma\xx;\tau)}\int_{V(P_\Gamma\xx;\tau)} \overline{K}(\xx,P_\Gamma\mathbf{z})\dd \sigma_{\mathbf{z}}.
\]
Then, for the kernels \eqref{eq:Kij-kernels}, the constants are
\begin{equation*}
  C^{(11)}_{\Gamma,\tau} = C^{(22)}_{\Gamma,\tau} = C^{(21)}_{\Gamma,\tau} = 0, \quad C^{(12)}_{\Gamma,\tau} = \dfrac{\exp(-\kappa \tau)-1+\kappa \tau}{2\pi \kappa \tau^2},
\end{equation*}
for $K_{11}$, $K_{22}$, $K_{21}$, and $K_{12}$ respectively.

    \subsection{The corrected trapezoidal rules in three dimensions - CTR}
    \label{subsec:CorrTR}
In this section we present how to construct the corrected trapezoidal rules used in \eqref{eq:corr-trapz-rule-sys1}-\eqref{eq:corr-trapz-rule-sys2}. 
The three dimensional quadrature rules will be defined
as the sum of integration over different coordinate planes. 
The two dimensional corrected trapezoidal rule 
is applied to approximate the integration over each plane. The rule belongs to a class of corrected trapezoidal rules in $\R^n$ for which convergence results have been proven \cite{izzo2022convergence}. The particular selection of the coordinate planes depends on the normal vector of $\Gamma$ at the target point located on the surface. 
To make the exposition clear, we will adopt the following convention for this section. 
\begin{notation}\label{not:sansserif-standard}
    In the standard basis $(\bar{\mathbf{e}}_1,\bar{\mathbf{e}}_2,\bar{\mathbf{e}}_3)$ of $\R^3$ we distinguish between variables in $\R^2$ and $\R^3$ by using sans serif variables for vectors in $\R^2$ and boldface variables with a bar for vectors of $\R^3$. For example, $\xs\in\R^2$ and $\bar\xx\in\R^3$.
    Moreover, given a point $\bar\xx$ in the standard basis, we denote by its sans serif character its first two components: 
    \[
    \bar\xx \equiv \left(\begin{array}{c} x_1 \\ x_2 \\ x_3 \end{array}\right) \equiv \left(\begin{array}{c} \xs \\ x_3 \end{array}\right).
    \]
\end{notation}

We consider a target point, $\bar\xx^*=(x^*,y^*,z^*)\in\Gamma$, at which
the surface normal is $\nn=(n_1,n_2, n_3)$, $|\nn|=1$. We define the dominant direction $\bar{\mathbf{e}}_i$ 
of the normal as the one with the largest component:
\[
i=\text{argmax}_{j=1,2,3}|n_j|,\quad (\bar{\mathbf{e}}_1,\bar{\mathbf{e}}_2,\bar{\mathbf{e}}_3)=\left(
\left(\begin{array}{c}
     1\\ 0 \\ 0 
\end{array}\right),\,
\left(\begin{array}{c}
     0\\ 1 \\ 0 
\end{array}\right),\,
\left(\begin{array}{c}
     0\\ 0 \\ 1 
\end{array}\right)
\right).
\]
Let $\rx$ be the unitary matrix describing the change of coordinates such that $\bar{\mathbf{e}}_i$ coincides with the $z$-direction. The matrix $\rx$ corresponds to a column permutation of the $3\times 3$ identity matrix. Let $\ry_3$ be the third row of $\rx$.
\begin{notation}\label{not:sansserif-newcoord}
    We denote the points in the new coordinate system obtained by the multiplication by $\rx$ using a tilde: given $\bar\xx$ expressed in the standard basis, 
    \[
        \tilde\xx := \rx \bar\xx
    \]
    is its expression in the new basis. Analogously to Notation~\ref{not:sansserif-standard}, given a point $\tilde \ux$ in the new coordinate system, we denote by its sans serif character its first two components: 
    \[
    \tilde\ux \equiv \left(\begin{array}{c} \tilde u_1 \\ \tilde u_2 \\ \tilde u_3 \end{array}\right) \equiv \left(\begin{array}{c} \uy{u} \\ \tilde u_3 \end{array}\right).
    \]
\end{notation}

For ease of notation, we define the projection mapping $\tilde P_\Gamma$ after change of coordinates:
\[
\tilde P_\Gamma(\tilde\xx) := P_\Gamma(\rx^{-1}\tilde\xx) = P_\Gamma(\bar\xx).
\]

For simplicity, we denote the integrands in \eqref{eq:PB-BIE-1}-\eqref{eq:PB-BIE-2}
using the function $f$,
\begin{equation}\label{eq:kernelfunction}
    f(\bar\y):={K}(\bar\xx^*, P_\Gamma(\bar\y))\rho(P_\Gamma(\bar\y))\delta_{\Gamma,\varepsilon}(\bar\y),\qquad
    \bar\y=(x,y,z).
\end{equation}
The dependence on $\bar\xx^*$ is suppressed in this notation, but we shall see how its role fits in the following derivation.

To approximate \eqref{eq:PB-BIE-1}-\eqref{eq:PB-BIE-2}
the standard trapezoidal rule is first applied in the dominant direction, which after change of coordinates corresponds to the third component.
With the grid points
$t_k=kh$, we get 
\begin{equation}\label{eq:plane-by-plane-sum}
    \int_{\R^3}f(\bar\y)\dd\bar\y  
    \approx 
    h\sum_{k}\int_{\R^{2}}f\left( \rx^{-1}\left(\begin{array}{c} \tilde x \\ \tilde y \\ kh \end{array}\right) \right)\dd \tilde x\dd \tilde y.
\end{equation}
We note that $f$ is singular in the new coordinate system along the line
\begin{equation}\label{eq:y0(z)_expression}
    \tilde\ux_0(t) = \rx\bar\xx^* +\dfrac{t-\ry_3\bar\xx^*}{\ry_3\nn}\rx\nn, \quad t\in \R.
\end{equation}
since 
\begin{equation*}
    \tilde P_\Gamma(\tilde\ux_0(t)) = P_\Gamma\left(\rx^{-1}\tilde\ux_0(t)\right) = P_\Gamma\left( \bar\xx^* +\dfrac{t-\ry_3\bar\xx^*}{\ry_3\nn}\nn\right) = P_\Gamma(\bar\xx^*), \quad \forall t\in\R.
\end{equation*}
By the assumption on $\rx\nn$, i.e. that the third component of $\rx\nn$ is the dominant, the normal does not lie on the plane $xy$. Then for any fixed $t$
\[
f\left( \rx^{-1}\left(\begin{array}{c} \cdot \\ \cdot \\ t \end{array}\right) \right)
\]
is singular at one point. 

In the following, for any fixed $t$ we derive the
the first two terms of a series expansion for each of 
the kernels \eqref{eq:Kij-kernels} around
the singularity.
Then a corrected trapezoidal rules for two dimensional integrals in \eqref{eq:plane-by-plane-sum} are derived based on the expansions.

From the definition of $\tilde\ux_0(t)$, its third component is $t$, and we define the points $\tilde\ux(t)$ as having third component equal to $t$:
\[
\tilde\ux_0(t)= \left(\begin{array}{c} \uy{u}_0(t)\\ t \end{array}\right), \text{ and}\quad \tilde\ux (t):=\left(\begin{array}{c} \uy{u}\\ t \end{array}\right), 
\]
where $\uy{u}_0(t)$ denotes the vector containing the first two components of $\tilde\ux_0(t)$, and $\uy{u}$ denotes the vector recording the first two coordinates of $\tilde\ux(t)$.
Then we factorize $f$, for a fixed $t$, as
\begin{equation}\label{eq:f_s_v_3D_splitting}
f\left(\rx^{-1} \tilde\ux (t) \right)= 
s(\uy{u}-\uy{u}_0(t);t)\, v\left(\tilde\ux (t)\right)
\end{equation}
where 
\begin{equation}\label{eq:svdefine}
\begin{array}{rl}
s(\uy{u};t)&={K}\left(\bar\xx^*, \tilde P_\Gamma\left(\begin{array}{c} \uy{u}+\uy{u}_0(t)\\ t \end{array}\right) \right),\\[0.3cm]
v\left(\begin{array}{c} \uy{u}\\ t \end{array}\right) &=\rho\left( \tilde P_\Gamma\left( \begin{array}{c} \uy{u}\\ t \end{array}\right)\right)\delta_{\Gamma,\varepsilon}\left( \rx^{-1}\left(\begin{array}{c} \uy{u}\\ t \end{array}\right)\right).
\end{array}
\end{equation}
Note that the type of singularity for $s$
depends on the properties of $\Gamma$ at the target point (such as principal curvatures, principal directions, normal).
Moreover, $s$ depends smoothly on $t$. 

We then use the second order corrected trapezoidal rule $\mathcal{U}_h^2$ to compute the integrals on each plane:
\begin{align*}
\int_{\mathbb{R}^2} f\left(\rx^{-1}\left(\begin{array}{c} \uy{u}\\ t_k \end{array}\right)\right)\dd\uy{u} &=
\int_{\mathbb{R}^2}
s(\uy{u}-\uy{u}_0(t_k);t_k) v\left(\begin{array}{c} \uy{u}\\ t_k \end{array}\right) \dd \uy{u}\\
&\approx \mathcal{U}^2_h \left[
s(\,\cdot-\uy{u}_0(t_k);t_k) v\left(\begin{array}{c} \cdot \\ t_k \end{array}\right) \right].
\end{align*}
We denote by ${\uy{u}}_{h}(t)$ and $(\alpha(t),\beta(t))$ 
the closest grid node to $\uy{u}_0(t)$ and the relative grid shift parameters
respectively, i.e.
\[
\uy{u}_0(t) = \uy{u}_h(t) + h\left(\begin{array}{c} \alpha(t)\\ \beta(t)\end{array} \right), \quad \alpha(t),\beta(t) \in \left[ -\frac12,\frac12 \right).
\]

From the definition in \eqref{eq:svdefine} we can compute the expansion 
\begin{equation}\label{eq:expansion-s-3D}
    s(\uy{u};t) = \dfrac{1}{|\uy{u}|}s_{0}\left(\dfrac{\uy{u}}{|\uy{u}|};t\right) + s_{1}\left(\dfrac{\uy{u}}{|\uy{u}|};t\right) + \Oo(|\uy{u}|).
\end{equation}
The relevant expressions for $s_k$ are given in 
Theorem~\ref{thm:kernelexpansions} below.
The second order rule is:
\begin{align}\label{eq:U2-ibim}
    \mathcal{U}^2_h \left[
    f\left(\rx^{-1}\left(\begin{array}{c} \cdot \\ t \end{array}\right)\right) \right] =&\, h^2\sum_{\uy{u}\in h\mathbb{Z}^2 \setminus \{\uy{u}_h(t)\}} f\left(\rx^{-1}\left(\begin{array}{c} \uy{u}\\ t \end{array}\right)\right) \\ & + h\, \omega[|\cdot|^{-1}s_{0}(\,\cdot\,; t);\alpha(t),\beta(t)]\,v\left(\begin{array}{c} \uy{u}_h(t) \\ t \end{array}\right).\nonumber
\end{align}
We define
\[
\bar\y_{h,k} := \rx^{-1}\left(\begin{array}{c} \uy{u}_h(t_k) \\ t_k \end{array}\right) = \rx^{-1} \tilde\ux_h(t_k) ,
\]
so the three-dimensional second order method $\mathcal{V}^{2,z}_h$, obtained by applying $\mathcal{U}^2_h$ \eqref{eq:U2-ibim} plane-by-plane along the dominant direction, is given by
\begin{align}
     \mathcal{V}^{2}_h[f] :=& h\sum_{k\in\mathbb{Z}} \mathcal{U}^2_h \left[
    f\left(\rx^{-1}\left(\begin{array}{c} \cdot \\ t_k \end{array}\right)\right) \right]  \label{eq:V2-3D-general} \\
     =& h^3\sum_{\bar\y\in h\mathbb{Z}^3\, \setminus\, \left(\bigcup_{k\in\mathbb{Z}} \{\rx^{-1} \tilde\ux_h(t_k)\}\right)}f(\bar\y) + h^2\sum_{k\in\mathbb{Z}} \omega[|\cdot|^{-1}s_{0}(\,\cdot\,;t_k);\alpha(t_k),\beta(t_k)]\,v\left(\begin{array}{c} \uy{u}_h(t_k) \\ t_k \end{array}\right)\nonumber \\
    =& h^3\sum_{\bar\y\in h\mathbb{Z}^3\, \setminus\, \left(\bigcup_{k\in\mathbb{Z}} \{\bar\y_{h,k}\}\right)}f(\bar\y) \nonumber + h^2\sum_{k\in\mathbb{Z}} \omega[|\cdot|^{-1}s_{0}(\,\cdot\,;t_k);\alpha(t_k),\beta(t_k)]
    \rho\left( \bar\y_{h,k} \right)\delta_{\Gamma,\varepsilon}\left( \bar\y_{h,k} \right).\nonumber 
\end{align}

        \subsubsection{Expansions of the Poisson-Boltzmann kernels}
        \label{subsec:expansions}

In this section we will analyze and expand as \eqref{eq:expansion-s-3D} the singular functions defined in \eqref{eq:svdefine} when $K$ are the Poisson-Boltzmann kernels
\eqref{eq:layer-kernels}. Specifically we will find the first expansion term $s_0$ for all kernels and in addition find $s_1$ for the kernel $K^{12}$ which has a different asymptotic behavior.
\begin{equation}\label{eq:layer-kernels-convenient-setting}
    \begin{array}{rl}
        s^{(11)}(\uy{u};t) &= \left\{\dfrac{\partial G_0}{\partial \n_y}-\dfrac{\varepsilon_E}{\varepsilon_I}\dfrac{\partial G_\ka}{\partial \n_y}\right\} 
        \left(\bar\xx^*,\tilde P_\Gamma \left( \begin{array}{c} \uy{u}+\uy{u}_0(t)\\ t \end{array} \right) \right)
        ,\\[0.4cm]
        s^{(12)}(\uy{u};t) &= \left\{ G_0-G_\ka\right\} \left(\bar\xx^*,\tilde P_\Gamma \left( \begin{array}{c} \uy{u}+\uy{u}_0(t)\\ t \end{array} \right) \right),\\[0.4cm]
        s^{(21)}(\uy{u};t) &= \left\{ \dfrac{\partial^2 G_0}{\partial \n_x \partial \n_y}-\dfrac{\partial^2 G_\ka}{\partial \n_x \partial \n_y}\right\}\left(\bar\xx^*,\tilde P_\Gamma \left( \begin{array}{c} \uy{u}+\uy{u}_0(t)\\ t \end{array} \right) \right),\\[0.4cm]
        s^{(22)}(\uy{u};t) &= \left\{ \dfrac{\partial G_0}{\partial \n_x}-\dfrac{\varepsilon_I}{\varepsilon_E} \dfrac{\partial G_\ka}{\partial \n_x}\right\} \left(\bar\xx^*,\tilde P_\Gamma \left( \begin{array}{c} \uy{u}+\uy{u}_0(t)\\ t \end{array} \right) \right).
    \end{array}
\end{equation}

We base the expansion on the local properties of $\Gamma$ around the target point $\bar\xx^*$. 
Let $\bar{\bm{\tau}}_1$, $\bar{\bm{\tau}}_2$, and $\nn$ be the principal directions and outward normal in $\bar\xx^*$. Let $M$ be the diagonal matrix of the principal curvatures, and after change of coordinates let ${A} \in\R^{2\times 2}$, $\mathbf{c},\mathbf{d}\in\R^{2\times 1}$, and $\alpha\in\R$ be the submatrices of the orthogonal change of basis matrix:
\begin{equation}\label{eq:def-Q}
    M:= \left(
        \begin{matrix}
        \ka_1 & 0 \\ 0 & \ka_2
        \end{matrix}
    \right), \qquad
    \left(
    \begin{array}{cc}
        {A} & \mathbf{c} \\ \mathbf{d}^T & \alpha
    \end{array}
    \right) := \left(
    \begin{array}{ccc}
        \text{\textemdash} & \rx\bar{\bm{\tau}}_1 & \text{\textemdash} \\ 
        \text{\textemdash} & \rx\bar{\bm{\tau}}_2 & \text{\textemdash} \\
        \text{\textemdash} & \rx\nn & \text{\textemdash}
    \end{array}
    \right).
\end{equation}
Let $\eta(t):=d_\Gamma(\rx^{-1}\tilde\ux_0(t))$, and 
\[
D_0(t) := (I-\eta(t)M)^{-1}.
\]
We now have all the definitions to express the expansion terms. 
\begin{theorem}\label{thm:kernelexpansions}
    Let $\bar\xx^*\in\Gamma$ be the target point.
    Suppose that 
    $\rx^{-1}\tilde\ux_0(t)\in T_\varepsilon$.
    Then, there is an $r>0$, depending on $t$, such that
    all the singular functions defined in \eqref{eq:layer-kernels-convenient-setting} can
    be written in the form
    \begin{equation}\label{eq:ellform}
        s^{(ij)}( \tilde{\mathsf{u}};t)=\dfrac{1}{|\tilde{\mathsf{u}}|}\ell^{(ij)}\left(| \tilde{\mathsf{u}}|,\dfrac{ \tilde{\mathsf{u}}}{| \tilde{\mathsf{u}}|};t\right), \quad |\tilde{\mathsf{u}}|<r,
        \qquad i,j=1,2,
    \end{equation}
    where $\ell^{(ij)}\in C^\infty((-r,r)\times {\mathbb S}^1)$.
    Moreover, the functions
    $s_0^{(ij)}(\tilde{\ys};t):\s^{1}\times\R\to\R$, $i,j=1,2$ and $s_1^{(12)}(\tilde{\ys};t):\s^{1}\times\R\to\R$
    in the expansion \eqref{eq:expansion-s-3D} are
    \begin{equation}\label{eq:s0s1-kernels} 
        \begin{array}{lll}
            & s_0^{(11)}(\tilde{\ys};t) &:= \left( 1-\dfrac{\varepsilon_E}{\varepsilon_I} \right)\dfrac{ \tilde{\ys}^T A^T {D_0(t)}^T M {D_0(t)} A \tilde{\ys} }{ 8\pi |{D_0(t)} A \tilde{\ys}|^3 } ,\\[0.3cm]
            & s_0^{(21)}(\tilde{\ys};t) &:= \dfrac{1}{8\pi|{D_0(t)} A \tilde{\ys}|},\\[0.3cm]
            & s_0^{(22)}(\tilde{\ys};t) &:= \left( 1-\dfrac{\varepsilon_I}{\varepsilon_E} \right)\dfrac{ \tilde{\ys}^T A^T {D_0(t)}^T M {D_0(t)} A \tilde{\ys} }{ 8\pi|{D_0(t)} A \tilde{\ys}|^3 },\\[0.5cm]
            & s_0^{(12)}(\tilde{\ys};t) &:= 0,\qquad s_1^{(12)}(\tilde{\ys};t) := \dfrac{\ka}{4\pi},
        \end{array} \qquad \text{with }\tilde{\ys}\in\s^{1},\ t\in\R.    
    \end{equation}
\end{theorem}
As previously mentioned, the kernel $K_{21}$ is the bottleneck for the order of accuracy. We can find more expansion terms for the function $s^{(ij)}$ in \eqref{eq:ellform}, and consequently increase the order of accuracy, if we improve the local approximation of the surface around the target point $\bar\xx^*$.

We rotate and translate $\Gamma$ so that $\bar\xx^*$ is in the origin and the principal directions and the normal coincide with the standard basis in $\R^3$. Then $\Gamma$ can be represented as the graph of a function $F$. In our current approach we are using the (pure) second derivatives of the function $F$, i.e. the principal curvatures. If we also use the third derivatives of $F$, we can write the expansion terms
\[
    s^{(11)}_1(\tilde{\ys};t), \quad s^{(22)}_1(\tilde{\ys};t), \quad s^{(12)}_2(\tilde{\ys};t),
\]
which we can use to build a third order accurate corrected trapezoidal rule, as we showed in \cite{izzo2022high}.
However, we cannot write one corresponding to $K_{21}$, i.e. $s^{(21)}_1(\tilde{\ys};t)$. To find such term, we need not only the third derivatives of $F$, but also the fourth derivatives.

        \subsubsection{Approximation and tabulation of the weights}
        \label{sub:weights}
Given the functions \eqref{eq:s0s1-kernels} we want to compute the weights $\omega[s_{0}^{(ij)}(\,\cdot\,;t_k);\alpha(t_k),\beta(t_k)]$. The functions in the expansion \eqref{eq:expansion-s-3D} are of the kind
\[
 \dfrac{1}{|\ys|}\ell\left(\dfrac{\ys}{|\ys|}\right), \quad \ys\in\R^2,
\]
for $s_0^{(ij)}$ where $i,j=1,2$ or, in the case of $s_1^{(12)}$, simply of the kind $\ell\left(\dfrac{\ys}{|\ys|}\right)$.
Fixed $(\alpha,\beta)$, we write $\ell$ using its Fourier series:
\begin{align*}
\ys=&|\ys|\big(\cos(\psi(\ys)),\sin(\psi(\ys))\big), \quad \ell\left(\dfrac{\ys}{|\ys|}\right)= a_0+\sum_{k=1}^\infty \big( a_k\cos(k\psi(\ys))+b_k\sin(k\psi(\ys)) \big),
\end{align*}
where $\{a_j\}_{j=0}^\infty$ and $\{b_j\}_{j=1}^\infty$ are the Fourier coefficients
of $\ell$. By linearity of the weights with respect to $\ell$, we can write them as
\begin{align*}
    \omega[s_0^{(ij)};\alpha,\beta] \,=&\, a_0\,\omega\left[|\ys|^{-1};\alpha,\beta\right]
    + \sum_{k=1}^\infty \left(a_k\,\omega\left[|\ys|^{-1}\cos(k\psi(\ys));\alpha,\beta\right]+b_k\,\omega\left[|\ys|^{-1}\sin(k\psi(\ys));\alpha,\beta\right]\right).
\end{align*}
We can then approximate and tabulate the weights 
$\omega[s_0^{(ij)};\alpha,\beta]$
in the following way. We fix a stencil of parameters 
$\{(\alpha_m,\beta_n)\}_{m,n}$
around $(\alpha,\beta)$ 
and {basis functions
}
$\{c_{m,n}(\alpha,\beta)\}_{m,n}$
such that we can approximate a function $f:\R^2\to\R$ in $(\alpha,\beta)$ as
\[
f(\alpha,\beta)\approx \sum_{m,n}c_{m,n}(\alpha,\beta)\,f(\alpha_m,\beta_n).
\]
We let $N$ be
the number of Fourier modes used to approximate the weights. 
Then, given $\ell$, we first find the $2N+1$ coefficients 
$a_0,\{a_k,b_k\}_{k=1}^N$ by using the Fast Fourier Transform.
Then, 
\begin{align*}
\omega[s_0^{(ij)};\alpha_m,\beta_n] \,\approx\, a_0\,\omega\left[|\ys|^{-1};\alpha_m,\beta_n\right] + \sum_{k=1}^N \Bigg(&a_k\,\omega\left[|\ys|^{-1}\cos(k\psi(\ys));\alpha_m,\beta_n\right]\\
& + b_k\,\omega\left[|\ys|^{-1}\sin(k\psi(\ys));\alpha_m,\beta_n\right]\Bigg),
\end{align*}
and we can approximate the weight for $(\alpha,\beta)$ via
\[
\omega[s_0^{(ij)};\alpha,\beta] \approx \sum_{m,n}c_{m,n}(\alpha,\beta)\,\omega[s_0^{(ij)};\alpha_m,\beta_n].
\]
We consequently need to compute and store the weights for the 
following constant and trigonometric functions,
\[
\left.
\begin{array}{l}
\omega\left[|\ys|^{-1};\alpha_m,\beta_n\right]\\[0.15cm]
\omega\left[|\ys|^{-1}\cos(k\psi(\ys));\alpha_m,\beta_n\right] \\[0.15cm] 
\omega\left[|\ys|^{-1}\sin(k\psi(\ys));\alpha_m,\beta_n\right]
\end{array}
\right\} \
\begin{array}{l}
    k=1,\dots,N,\\
    \text{and all $m,n$ in the stencil for $(\alpha,\beta)$}.
\end{array}
\]
For $s_1^{(12)}$, the function is constant with respect to $|\ys|$ and $\ys/|\ys|$, so the only weights needed are $\omega\left[1;\alpha_m,\beta_n\right]$.
For ease of notation, let's write $s_q(\ys):=|\ys|^{q-1}\ell(\ys/|\ys|)$, where $q\in\{0,1\}$.
The weights $\omega[s_q;\alpha_m,\beta_n]$ are formally the limits 
\begin{equation}
  \omega[s_q;\alpha,\beta] := \lim_{h\to0^+} \omega_h[s_q;\alpha,\beta]\,,
  \label{eq:single_correction_weight}
\end{equation}
where 
\begin{align*}
  \omega_h[s_q;\alpha,\beta] :=& 
  \dfrac{1}{h^{q+1}}\frac{ \int_{\mathbb{R}^2}{s_q}(\xx){g}(\xx)\text{d}\xx-T_{h}^{0}\Big[s_q(\,\cdot -(\alpha,\beta)h) g(\,\cdot -(\alpha,\beta)h)\Big]}{ g( - (\alpha,\beta)h)}.
\end{align*}
The function $g$ is chosen to be $g \in C^\infty_c({\mathbb R}^2)$ and radially symmetric, with $g(\zero)=1$. 

\begin{remark}
    We approximate the limit \eqref{eq:single_correction_weight} by $\omega_{i,h^*}$, where
    \[
    h^*:=2^{-n}, \ \ n := \arg \min_{j=1,2,3,\dots} \left\{ |\omega_{i,2^{-j}}-\omega_{i,2^{-j-1}}|\leq \texttt{Tol}\right\}.
    \]
    In the simulations presented in Section~\ref{sec:num:smooth}, to compute $\omega$, we use $\texttt{Tol}=10^{-8}$.
\end{remark}

    \subsection{Overview of the algorithm}
    \label{subsec:overview-algorithm}

\begin{figure}[!htb]
    \centering
    \includegraphics[scale=1]{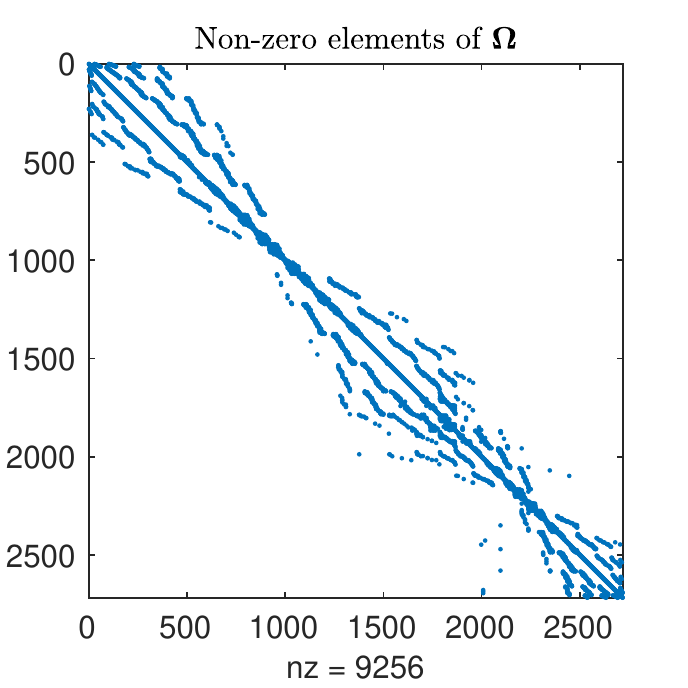}\includegraphics[scale=1]{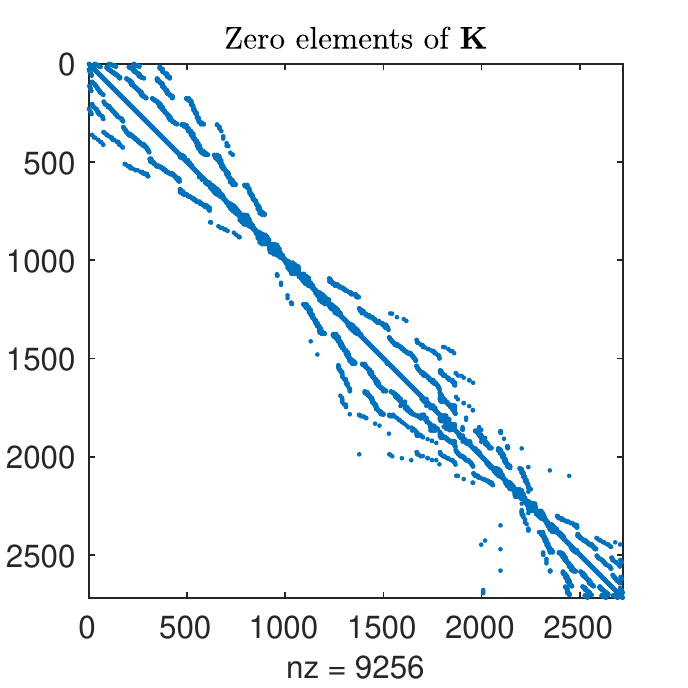}
    \caption{ Sparsity pattern of the matrix $\mathbf{\Omega}$ and complementary form of the matrix $\mathbf{K}$. 
    The surface is the SES (see Section~\ref{subsec:Surface-Definition}) for the union of two intersecting spheres of different radii. The tubular neighborhood is of width $\varepsilon=2h$, which means that in the best possible case 5 nodes are corrected for a given target node. 
    The weights present in $\mathbf{\Omega}$ substitute the removed values from $\mathbf{K}$ and correct them according to the CTR2 rule. 
    Left: non-zero elements of $\mathbf{\Omega}$, which has approximately 4 non-zero elements for each row. 
    Right: the zero elements of the dense matrix $\mathbf{K}$ corresponding to the ones removed in order to apply the correction. }
    \label{fig:spy-matrices}
\end{figure}

We now have all the tools to build the system \eqref{eq:corr-trapz-rule-sys1}-\eqref{eq:corr-trapz-rule-sys2}. Let $\{\y_k\}_{k=1}^M$ be the discretization nodes inside the tubular neighborhood. We write the system in matrix form as
\begin{equation}\label{eq:linear-system}
    \Lambda \mathbf{p} + h^3\mathbf{K}\mathbf{W}\mathbf{p} + h^2\mathbf{\Omega}\mathbf{W}\mathbf{p} = \mathbf{g}
\end{equation}
where $\mathbf{p}$ contains both $\bar\psi(\y_k)$ and $\bar\psi_n(\y_k)$, and
\[
\Lambda := 
\left(
\begin{array}{cc}
    \lambda_1 \mathbf{I} & \mathbf{0} \\ \mathbf{0} & \lambda_2 \mathbf{I} 
\end{array}
\right).
\]
The matrix $\mathbf{W}$ is a diagonal matrix defined by the weights $\delta_{\Gamma,\varepsilon}(\y_k)$,
and the matrices $\mathbf{K}$ and $\mathbf{\Omega}$ represent the dense matrix of the kernel evaluations $\overline{K}(P_\Gamma\y_k,\y_m)$, $k,m=1,\dots,M$ except in the corrected nodes, and the sparse matrix of the weights for the corrected nodes respectively. In Figure~\ref{fig:spy-matrices} we see the complementary structure of the matrices $\mathbf{K}$ and $\mathbf{\Omega}$, where the missing elements in $\mathbf{K}$ correspond to the corrected terms in $\mathbf{\Omega}$.  

We solve this system by a standard \texttt{GMRES} solver. In the \texttt{GMRES} algorithm, we use the black-box fast multipole method \texttt{BBFMM} (see \cite{fong2009black}) to accelerate the multiplication of the operator $\mathbf{K}$ to any vector. 
The codes are available on GitHUB\footnote{\texttt{https://github.com/lowrank/ibim-levelset}}.

Compared to the K-reg method, the {CTR2} has a similar computational cost. 
The weights require interpolation of the tabulated values, which is negligible. For the second order method described here, CTR2, only one node per plane needs correction, which is similar to K-reg is the parameter $\tau=h$ or $2h$. 
The corrected nodes for CTR2 and the regularized nodes for K-reg lead to a matrix $\Omega$ in \eqref{eq:linear-system} which has a nonzero diagonal and at most $2\varepsilon/h$ other nonzero elements for each row. 
By taking $\varepsilon=Ch$, the number of nonzero row elements is independent of $h$.

The {CTR2} needs to identify all the closest nodes along the singular line and interpolate the weights. 
This is done once in precomputation, as the matrix used in the matrix-vector multiplication does not change along the \texttt{GMRES} iterations.

In the following Section, we will present the number of \texttt{GMRES} iterations needed to solve the system for the methods used, and see that {CTR2} has a similar condition number compared to K-reg, needing a small number of iterations to converge to the desired accuracy.

    \subsection{Numerical tests}
    
\label{sec:num:smooth}

In this section we present two numerical studies on the order of accuracy of our method (CTR2). We compare the results to those computed with the IBIM regularization (K-reg) method
described in Section~\ref{subsec:Kreg}. 

The surfaces used in our tests include a sphere with radius $r=10$ and a torus with radii $R_1=1$ and $R_2=1/2$, both centered in the origin. The torus is rotated along the $x$-, $y$-,and $z$-axes by the following angles respectively:
\[
    a = 1.99487,\quad b = 2.54097947651017, \quad c = 4.219760487439292.
\]
We assess our method by 
solving numerically \eqref{eq:PB-BIE-1}-\eqref{eq:PB-BIE-2} on the grid $h\mathbb{Z}^3$ with right hand side $g$ built using a single charge at $\zero$ with charge $q$. 
The constants $\varepsilon_I$ and $\varepsilon_E$ are fixed at:
\[
    \varepsilon_I = 1.0, \quad \varepsilon_E = 80.0.
\]
Then we find the solutions 
$\psi$ and $\psi_n=\partial \psi/\partial \n$.
Then we compute the integral 
\begin{align}
  \psi_{rxn}(\mathbf{z}) &:= \int_\Gamma \left\{ \left( \frac{\varepsilon_E}{\varepsilon_I} \dfrac{\partial G_\kappa}{\partial \n_y}(\mathbf{z},\y)-\dfrac{\partial G_0}{\partial \n_y}(\mathbf{z},\y) \right)\psi(\y) + \left( G_0(\mathbf{z},\y)-G_\kappa(\mathbf{z},\y) \right)\dfrac{\partial \psi}{\partial \n}(\y) \right\}\dd\sigma_\y, \label{eq:psi-rxn}
\end{align}
where $\mathbf{z}$ is a point in space not belonging to the surface. The integral is computed in the IBIM framework \eqref{eq:surface-integral-IBIM} using the standard trapezoidal rule for
\[
f(\y) := \left( \frac{\varepsilon_E}{\varepsilon_I} \dfrac{\partial G_\kappa}{\partial \n_y}(\mathbf{z},\y)-\dfrac{\partial G_0}{\partial \n_y}(\mathbf{z},\y) \right)\psi(\y) + \left( G_0(\mathbf{z},\y)-G_\kappa(\mathbf{z},\y) \right)\dfrac{\partial \psi}{\partial \n}(\y),
\]
and the averaging kernel $\phi$ \eqref{eq:phi-averaging-kernel} used is
\[
\phi(x) = 
\begin{cases}
\dfrac{1}{2}(1+\cos(\pi x)), & \text{if }|x|\leq 1, \\[0.2cm]
0, & \text{otherwise}.
\end{cases}
\]
We use the standard trapezoidal rule because $\mathbf{z}$ does not lie close to the surface, so the integrand $f$ is smooth.

\begin{figure}[!htb]
    \centering
    \includegraphics[scale=1]{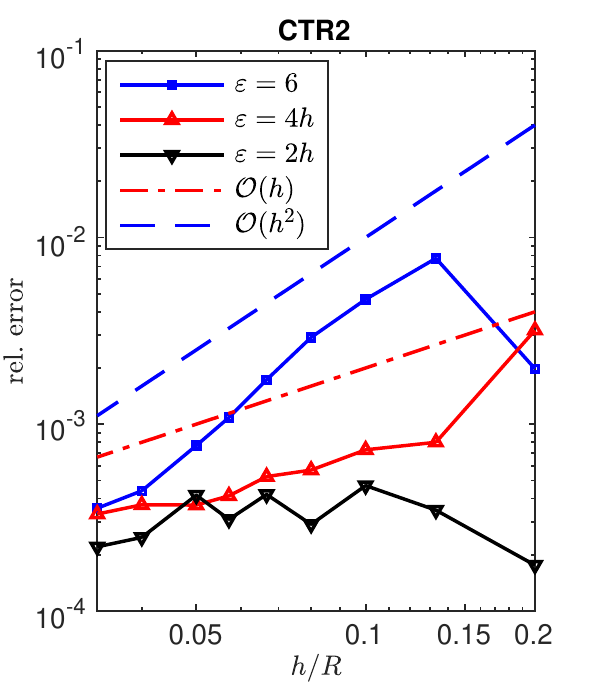}\includegraphics[scale=1]{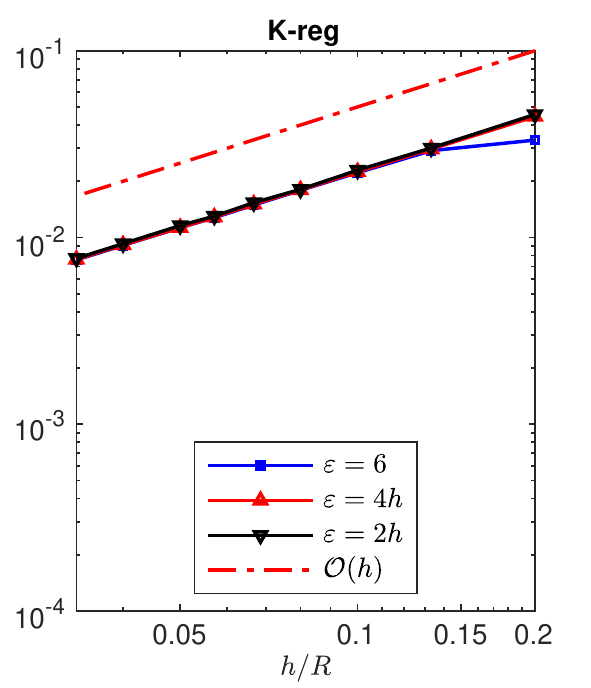}
    \caption{Convergence results for a sphere. The methods K-reg and CTR2 of order two are used to discretize the integrals in \eqref{eq:PB-BIE-1}-\eqref{eq:PB-BIE-2} and solve the corresponding system. Then we evaluate $\psi_{rxn}(\xx_0)$ in $\xx_0$ center of the sphere. For the sphere we have analytic values for $\psi_{rxn}(\xx_0)$. We plot the relative error for CTR2 (left) and K-reg (right).}
    \label{fig:single-sphere-ion}
\end{figure}

In Figure~\ref{fig:single-sphere-ion} we compare the convergence of the error in $\psi_{rxn}(\zero)$ for the sphere, for K-reg and CTR2 and different widths of the tubular neighborhood $T_\varepsilon$.
The exact values $\psi^*$ and $\psi_n^*$ of the solution on the sphere are:
\begin{align*}
    \psi^*(\xx) &= 
    \begin{cases}
    \dfrac{q}{4\pi\varepsilon_E |\xx|}+\dfrac{q}{4\pi r}\left( \dfrac{1}{\varepsilon_E(1+\kappa r)}-\dfrac{1}{\varepsilon_I} \right), & |\xx| \leq r, \\
    \dfrac{q e^{-\kappa(|\xx|-r)}}{4\pi\varepsilon_E (1+\kappa r)|\xx|}, & |\xx| > r.
    \end{cases}\\
    \psi^*_n(\xx) &= 
    \begin{cases}
    -\dfrac{q}{4\pi\varepsilon_I |\xx|^2}, & |\xx| \leq r, \\
    -\dfrac{q e^{-\kappa(|\xx|-r)}}{4\pi\varepsilon_E (1+\kappa r)|\xx|^2}-\dfrac{\kappa q e^{-\kappa(|\xx|-r)}}{4\pi \varepsilon_E (1+\kappa r)|\xx|}, & |\xx| > r,
    \end{cases}
\end{align*}
and the value of $\psi_{rxn}(\zero)$ is
\[
\dfrac{q^2}{4\pi r}\left( \dfrac{1}{\varepsilon_E(1+\kappa r)}-\dfrac{1}{\varepsilon_I} \right).
\]

The regularization parameter used for K-reg is $\tau=2h$. 
When $\varepsilon$ is constant with respect to $h$, CTR2 exhibits an order of accuracy slightly greater than 2, while K-reg only has order 1. 
In the case $\varepsilon=2h$, the order of accuracy of CTR2 decreases to approximately 1 
{because $\delta_{\Gamma,\epsilon}$is not well resolved by the grid.  }
Then the order of accuracy is the same as K-reg; however the error constant for CTR2 is smaller by a factor of at least 10 compared to K-reg.\\

\begin{table}[!htb]
    \caption{{Torus test, with $\varepsilon=2h$.} We compute the potentials $\psi$ and $\psi_n$ numerically using K-reg and CTR2, and use them to compute  $\psi_{rxn}(\zero)$ \eqref{eq:psi-rxn}. We use two different approximations of the Jacobian $J_\Gamma$. For K-reg, $J_\Gamma\approx 1$, and for CTR2, $J_\Gamma$ is approximated to second order in $h$ by $J_\Gamma\approx J_h$. 
    For the surface area $\mathcal{A}$ we directly tabulate the relative error. }
    \begin{center}
        \begin{tabular}{c|c|c|c|c|c|c}
            grid size & {$\text{err}_\mathcal{A}$}, $J_\Gamma \approx 1$ & $\psi_{rxn}(\zero)$, K-reg & diff & {$\text{err}_\mathcal{A}$}, $J_\Gamma \approx J_h$ & $\psi_{rxn}(\zero)$, CTR2 & diff \\ \hline
            103 & 9.483e-6 & 2.213e+2 &        & 8.521e-6 & 2.2359e+2 &           \\
            128 & 9.280e-6 & 2.219e+2 & 5.6e-1 & 4.498e-6 & 2.2354e+2 & 4.9e-2 \\
            256 & 2.283e-5 & 2.229e+2 & 9.8e-1 & 9.525e-7 & 2.2353e+2 & 1.3e-2  \\
            512 & 1.783e-6 & 2.232e+2 & 3.3e-1 & 2.492e-6 & 2.2352e+2 & 2.5e-3 
        \end{tabular}
    \end{center}
    \label{tab:torus}
\end{table}

For the torus we compiled the computed values in Table~\ref{tab:torus}. 
In the computation of $\psi_{rxn}(\zero)$, K-reg approximates the Jacobian in the IBIM integral \eqref{eq:surface-integral-IBIM} as constant 1, while CTR2 approximates it using 
a second order approximation in $h$, $J_\Gamma\approx J_h$. 
{According to \cite{kublik2018extrapolative}, if $\phi$ in \eqref{eq:phi-averaging-kernel} is an even function, then the terms corresponding to the odd powers of the distance in the Jacobian, $J_\Gamma$, will average out analytically. 
Hence, by approximating $J_\Gamma$ by $1$, one may expect an analytically error of $\mathcal{O}(\varepsilon^2)$. 
This fact was exploited for convenience in K-reg.
However, when $\varepsilon$ is only a small constant multiple of $h$, discrete errors will dominate (see \cite{engquist2005discretization}). 
In CTR2, since curvature information is required for computing the quadrature weights, it is natural to use a second order approximation of $J_\Gamma$.}
To gauge the influence of these approximations, we also compute the surface area $\mathcal{A}$ by using formula \eqref{eq:surface-integral-IBIM} with $f\equiv 1$ and the same two ways of approximating $J_\Gamma$. 
We then compute the relative surface area error, listed in Table~\ref{tab:torus}.

\section{Computing electrostatic potential of macromolecules}
\label{sec:num:application}

In this section, we use the IBIM formulation to solve the linearized Poisson-Boltzmann equation and compute the electrostatic potential and polarization energy for solvent-molecule interfaces. 
The solvent-molecule interface represents the interface separating the solvent fluid particles and the macromolecule.
The boundary integral equations are defined on such interface.
Such surfaces are naturally complex and difficult to parametrize. Consequently it is difficult to apply standard BIM with an explicit parametrization for these applications.
In the following section we briefly review two different approaches for defining the solvent-molecule interface. For either approach the IBIM approach can be conveniently applied.
Then the IBIM equations are discretized using CTR. We compare the method with the previous regularization approach K-reg, and introduce a hybrid corrected/regularized method HYB-IBIM to deal with surfaces which are piecewise smooth.

    \subsection{{Solvent-molecule interface }}
    \label{subsec:Surface-Definition}
The solvent-molecule interface, $\Gamma$, can be defined in different ways. 
Classically, $\Gamma$ is approximately by unions of predefined shapes, which leads to 
less accurate approximations.
Nevertheless, one can obtain qualitatively correct information when the errors are averaged over large molecules \cite[\S 22.1.2.1.2]{gerstein2006protein}. 

The van der Waals (VDW) molecular surface is the simplest interface to define: each atom in the protein is approximated by a solid sphere with the van der Waals radius. 
{The union of the regions enclosed by the spheres is the VDW surface.} 
A VDW surface typically has kinks are is only $C^0$. 

{The solvent-accessible surface (SAS) of a molecule is defined by "rolling" a small sphere around the VDW surface. A small region near the kink in the VDW surface are replaced by a portion of this sphere, and the the resulting SAS is $C^1$ and piece-wise $C^2$.}
From the SAS one can define the solvent-exluding surface (SES) by removing any cavity enclosed by the molecule, and consequently inaccessible by the solvent. The radius of the probe is usually set to the radius of the water molecule. The accessibility of solvent actually depends on its three-dimensional orientation, but this dependence is ignored for simplicity because the radius of the oxygen atom is much larger than the hydrogen 
(see \cite[\S 22.1.2.1.2]{gerstein2006protein}). The resulting SES is formally $C^1$.

The SES can be quite difficult to find analytically, so we used the numerical approach described in \cite{zhong2018implicit}. The information about the SES is in the form of a signed distance function on a uniform Cartesian grid. The steps are summarized here:
\begin{enumerate}
  \item An initial level set function is defined.
  \item The level set function is evolved using an ``inward'' eikonal flow. 
  \item The internal cavities are removed. 
  \item A reinitialization procedure is applied to the level set function, which will give as a result the sought signed distance function to the SES.
\end{enumerate}

The physical phenomena which control the dynamics of the {biomolecular} systems greatly impact the solvent-molecule interface. 
Consequently in biomolecular applications, {for improving accuracy to the underlying physics,} the solvent-molecule interface itself can be an unknown to solve for. 
More recently methods have been developed which take the relevant dynamics into account, and find the solvent-molecule interface as the minimizer of an energy functional. One such method is the Variational Implicit Solvation Model (VISM).

\begin{figure}[!htb]
    \centering
    \includegraphics{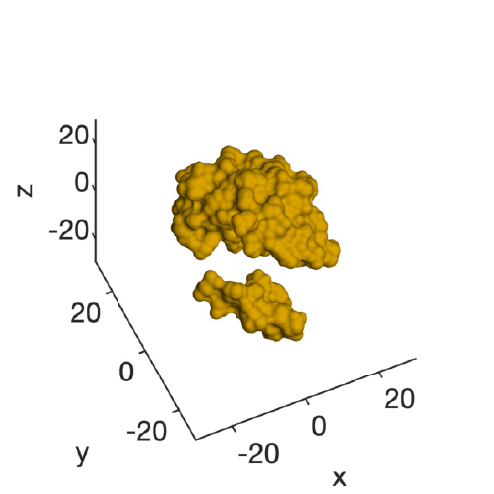}\includegraphics{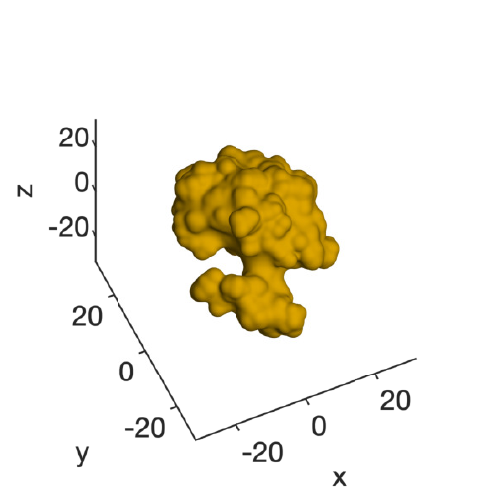}\includegraphics{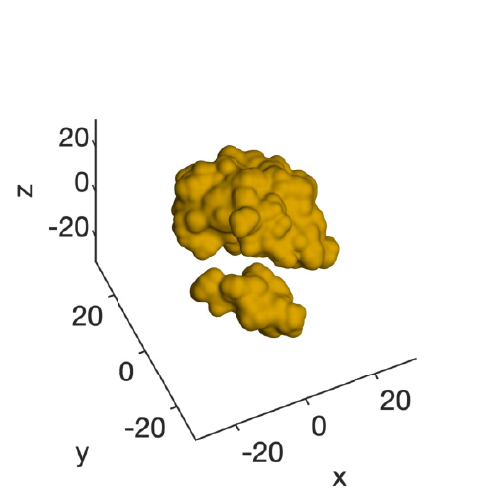}
    \caption{
    Solvent-molecule interfaces of p53-MDM2 generated in two different ways. On the left the SES is generated using the code from \cite{zhong2018implicit}. 
    The center and right figures show the smooth solvent-molecule interfaces generated via VISM. 
    The two different surfaces are generated using two different initial configurations of the complex biomolecular system p53-MDM2 (PDB ID 1YCR) \cite{kussie1996structure} from the Protein Data Bank (PDB), as detailed in \cite{zhang2021coupling}. 
    The different initial configurations can lead to different minimizers to the free-energy functional because of its non-convex nature. 
    This is detailed in the original algorithm paper of VISM \cite{wang2012level}.
    The surfaces are plotted using the signed distance function values. Left: SES surface. Center: VISM ``loose'' surface, generated with an initial configuration further from the atoms. Right: VISM ``tight'' surface, generated with an initial configuration closer to the atoms.}
    \label{fig:smooth-molecule-1}
\end{figure}

In VISM, the solvent-molecule interface is found by iteratively minimizing the free-energy functional over a surface, as for example seen in \cite{wang2012level}:
\begin{align*}
    \mathcal{G}[\Gamma] 
    &=  \mathcal{G}_{geom}[\Gamma] + \mathcal{G}_{vdW}[\Gamma] + \mathcal{G}_{elec}[\Gamma],
\end{align*}
where the electrostatic interaction $\mathcal{G}_{elec}[\Gamma]$ is approximated using the Coulomb-field approximation.
Typically, this energy is minimized by computing a sequence of minimizing surfaces, $\Gamma_k$, $k=0,1,2,\cdots.$ For each $k$ an electrostatic problem is solved {to evaluate $\mathcal{G}_{elec}[\Gamma_k]$}. While other formulations of the electrostatic energy can be used, e.g. the Coulomb-field approximation (CFA) (see e.g. \cite{wang2012level,zhang2021binary}) one of the most accurate is provided by the Poisson-Boltzmann theory (see e.g. \cite{li2011dielectric, zhou2014variational, li2015diffused}).

Different methods lead to solvent-molecule interfaces of different smoothness: 
\begin{itemize}
    \item The SAS is of smoothness $C^0$ or $C^1$ depending on which set of points is used. 
    In the case of the SAS defined as the set of center points of the probe as it rolls the surface can have kinks and consequently be of regularity at most $C^0$. 
    In the case of the SAS defined as the set of points tangent to the model the surface is of regularity $C^1$ because the curvatures will be discontinuous when jumping from an atom to the probe.
    \item The SES is formally of the same smoothness of the SAS taken with points tangent to the model.
    \item The surface generated with VISM is smoother due to the presence of the surface tension term in $\mathcal{G}_{geom}[\Gamma]$, one of the terms in the energy to be minimized.
\end{itemize}

{We point out that IBIM requires that the solvent-molecule interface to be at least $C^{1,\alpha}$ for some $\alpha>0$.
This way, one can be sure that there exists a workable tubular neighborhood. 
The second order CTR, CTR2, formally requires the surfaces to be $C^2$ almost everywhere on the surface, since principal curvatures and directions on the surface are needed for computing the quadrature weights. 
While SES surfaces formally do not present any problem on the analytical level,
the approximation of these second order geometrical quantities on an SES, near the discontinuities of these surface quantities, may introduce non-negligible errors to the integration and to the discretized linear system.}
In the following subsections, we propose a hybrid method for solving the linearized Poisson-Boltzmann equation on the SES. 
We will show that the hybrid method retains second order of accuracy if, for a piecewise smooth surface, the number of points where the curvatures are undefined scales linearly with $h$.

\subsection{Hybrid corrected/regularized quadrature}\label{sec: hyd cor quad}

We classify the nodes inside the tubular neighborhood
into nodes which are ``bad''  and ``good'', 
\[
    \mathcal{M} := \{\xx_m\}_{m=1}^M = h\mathbb{Z}^3 \cap T_\varepsilon = \mathcal{M}_B \cup \mathcal{M}_G,
\]
depending on whether the stencil to compute the second derivatives of the signed distance function is negatively impacted by the curvature discontinuity. 
By taking the corrected trapezoidal rule for the ``good'' points where the curvatures are relatively accurate, and the regularized quadrature rule for the ``bad'' points, we expect a second order accuracy in the grid size $h$ for the good points and only first order for the bad points. {The overall order of accuracy of the hybrid approach will then depend on the percentage of the bad points.}

For good and bad target points, the quadrature rules used to discretize the equations \eqref{eq:PB-BIE-1}-\eqref{eq:PB-BIE-2} are CTR2 and K-reg respectively. Then the system \eqref{eq:linear-system}, after reordering so that the nodes are grouped together, becomes
\begin{equation}\label{eq:linear-system-hyb}
    \left(
    \begin{array}{c}
        \Lambda_{|\mathcal{M}_{\mathcal{B}}|}+ h^3\mathbf{K}_{\mathcal{B}}\mathbf{W}_{\mathcal{B}} + h^2\mathbf{\Omega}_{Kreg}\mathbf{W}_{\mathcal{B}} \\
        \Lambda_{|\mathcal{M}_{\mathcal{G}}|}+ h^3\mathbf{K}_{\mathcal{G}}\mathbf{W}_{\mathcal{G}} + h^2\mathbf{\Omega}_{CTR}\mathbf{W}_{\mathcal{G}}
    \end{array}  
    \right)
    \left(
    \begin{array}{c} 
        \mathbf{p}_{\mathcal{B}}\\ 
        \mathbf{p}_{\mathcal{G}}
    \end{array}  
    \right) = 
    \left(
    \begin{array}{c}
        \mathbf{g}_{\mathcal{B}}\\ 
        \mathbf{g}_{\mathcal{G}}
    \end{array}  
    \right).
\end{equation}
The vector $\mathbf{p}_{\mathcal{B}}$ groups the values of $\psi$ and $\psi_n$ in the ``bad'' points, and analogously for $\mathbf{p}_{\mathcal{G}}$.
Depending on whether the target point is in $\mathcal{M}_{\mathcal{G}}$ or $\mathcal{M}_{\mathcal{B}}$, the nodes which are corrected or regularized may differ. Consequently we write the subscript $\mathcal{G}$ and $\mathcal{B}$ respectively.

The scheme will have order of accuracy $\Oo(h^{3/2})$ if the number of bad points scales at most linearly in $h^{-1}$. We see this by rewriting the system \eqref{eq:linear-system-hyb} in the simplified representation
\[
(\bm{\lambda} I + h^3 \mathbf{A}) \psi = \mathbf{b}
\]
where $\mathbf{A}$ represents the quadrature discretization matrix, $\psi$ represents the whole solution vector and $\mathbf{b}$ the right hand side. 

Given two nodes $\xx_i \in\mathcal{M}_G$ and $\xx_j \in \mathcal{M}_B$ in the good and bad sets respectively, 
then we can formally estimate the error of the quadrature rule as 
\begin{align*}
    h^3 \mathbf{A}(i,:)\psi &= \mathcal{S}[\psi](\xx_i) + \Oo(h^2), \\
    h^3 \mathbf{A}(j,:)\psi &= \mathcal{S}[\psi](\xx_j) + \Oo(h),
\end{align*}
where 
\[
    \mathcal{S}[\psi](\xx_i)=:\int_{T_\varepsilon} \overline{K}(P_\Gamma\xx_i,\y) \psi(\y) \delta_{\Gamma,\varepsilon}(\y) \dd \y 
\]
is the linear integral operator we are approximating. We call $\tilde \psi$ the solution to the system
\[
    \bm{\lambda} \tilde\psi_i + \mathcal{S}[\tilde \psi](\xx_i) + \mathbf{d}_i(h) = \mathbf{b}, \quad i=1,\dots,M,
\]
where
\[
    \mathbf{d}(h) := (\underbrace{h,h,\dots,h}_{\mathcal{M}_B},\underbrace{h^2,h^2,\dots,h^2}_{\mathcal{M}_G}),
\]
and $\psi$ is the solution to the system
\[
    \bm{\lambda} \psi_i + \mathcal{S}[ \psi](\xx_i) = \mathbf{b}, \quad i=1,\dots,M.
\]
Then the difference between the two solutions is:
\begin{align*}
    \psi-\tilde\psi = (\bm{\lambda} I+\mathcal{S})^{-1} \mathbf{d}(h).
\end{align*}
In the case $M\sim h^{-2}$ and $|\mathcal{M}_B|\sim h^{-1}$, then the dependency of the mean $\ell_2$-error of $\psi-\tilde\psi$ on $h$ is

\begin{align*}
    &{|\psi-\tilde\psi|} \leq |(\bm{\lambda} I+\mathcal{S})^{-1}| |\mathbf{D}(h)| = C \sqrt{\sum_{i=1}^{|\mathcal{M}_B|} h^2 + \sum_{i=1}^{|\mathcal{M}_G|} h^4} =  C \sqrt{ h^{-1} h^2 + h^{-2} h^4} \sim \Oo( \sqrt{h} )  \\
    \Longrightarrow & \quad \frac{|\psi-\tilde\psi|}{\sqrt{M}} \sim \Oo( h^{3/2} ).
\end{align*}
Then the averaged order of accuracy for the hybrid rule is formally $3/2$, still larger than $1$, if the number of bad points scales linearly with $h^{-1}$. \\

We call this hybrid method {HYB}-IBIM, and remark that the good and bad points are also used to identify which approximation of the Jacobian $J_\Gamma$. In the good points we use the second order approximation $J_\Gamma\approx J_h$, while in the bad points we use the approximation $J_\Gamma\approx 1$. We indicate this approximation of the Jacobian \eqref{eq:jacobian} by
\begin{equation}\label{eq:Jacobian-hybrid}
    J^{(H)}_h(\y) := \begin{cases}
    1, & \text{if }\y\in\mathcal{M}_B, \\
    J_h(\y), & \text{if }\y\in\mathcal{M}_G.
    \end{cases}
\end{equation}

    \subsection{Numerical tests}
    In the following, we assess our methods  by comparing the computed values of surface area and the polarization energy \eqref{eq:pol-energy}. 
The surface area, $\mathcal{A}$, is computed by {the standard trapezoidal rule} for the integral defined on the right-hand-side \eqref{eq:surface-integral-IBIM}  with $\rho$ constant equal to $1$ and $\phi$ as in Section~\ref{sec:num:smooth}. 
The polarization energy is defined by  the formula:
\begin{align}
  \mathcal{G}_{pol} &:= \frac12 \sum_{k=1}^{N_c} q_k \psi_{rxn}(\mathbf{z}_k), \label{eq:pol-energy} \\
  \psi_{rxn}(\mathbf{z}) &:= \int_\Gamma \left\{ \left( \frac{\varepsilon_E}{\varepsilon_I} \dfrac{\partial G_\kappa}{\partial \n_y}(\mathbf{z},\y)-\dfrac{\partial G_0}{\partial \n_y}(\mathbf{z},\y) \right)\psi(\y) + \left( G_0(\mathbf{z},\y)-G_\kappa(\mathbf{z},\y) \right)\dfrac{\partial \psi}{\partial \n}(\y) \right\}\dd\sigma_\y, \nonumber
\end{align}
where $\psi$ and $\partial \psi/\partial \n$ are computed by solving \eqref{eq:corr-trapz-rule-sys1} and \eqref{eq:corr-trapz-rule-sys2}. The surface integral \eqref{eq:pol-energy} is, again, computed by second order CTR, CTR2, applied to the associated IBIM formulation. Because $\mathbf{z}_k$ are the centers of the atoms, they do not fall onto the surface and the CTR2 just becomes the standard trapezoidal rule.

\subsubsection{Smooth interfaces}

In this section, we apply CTR2 to solve the linearized Poisson-Boltzmann equation on smooth surfaces. The tests run in Section~\ref{sec:num:smooth} include a sphere, which corresponds to a single ion. The value \eqref{eq:psi-rxn} corresponds to the polarization energy (up to a constant), and we presented convergence results in Figure~\ref{fig:single-sphere-ion}.

Then we apply CTR2 to solve systems derived from two smooth solvent-molecule interfaces for the protein p53-MDM2 \cite{zhang2021coupling}. The interfaces are constructed by the LS-VISM code developed in \cite{zhou2015ls} using two different starting configurations (corresponding to what is called a "loose" and a "tight interface).
Because of the complexity of the problem and the protein, we do not have analytical solutions.

The signed distance functions to the two surfaces
are computed on the Cartesian grid in $[-31.7945,31.7945]^3$, accurate within the absolute distance of $1.24197$ to each interface. 

We see in Figure~\ref{fig:smooth-molecule-1} that the tight surface has two closely positioned connected components.
In physical units we can estimate the bound on the reach to be approximately 2.5. 
{The stencil of the finite differences used to approximate the curvatures to second order in $h$ is 5 grid nodes in each coordinate direction. This means that we the $513^3$ uniform Cartesian grid ($h\approx 0.12$) can provide roughly 10 points to discretize the tight spacing between the two component, supporting a tubular neighborhood of width $\varepsilon=2h$. } 

In Tables~\ref{tab:smooth-surface-VISM-1} and~\ref{tab:smooth-surface-VISM-2} 
we tabulate the surface area and solve the system with right hand side computed using a single fictitious ion at the origin, 
which falls inside the loose molecule but outside the tight one. Then we evaluate $\psi_{rxn}$ \eqref{eq:psi-rxn} at this point, and its sign depends whether the evaluation point lies inside or outside of the surface.
The results are obtained for $\varepsilon=2h$ and $5h$, and we use the difference between them to get a rough estimate of the accuracy of the values.  

\begin{figure}
    \centering
    \includegraphics[trim = 3cm 0cm 0cm 1.8cm, clip, scale=0.4]{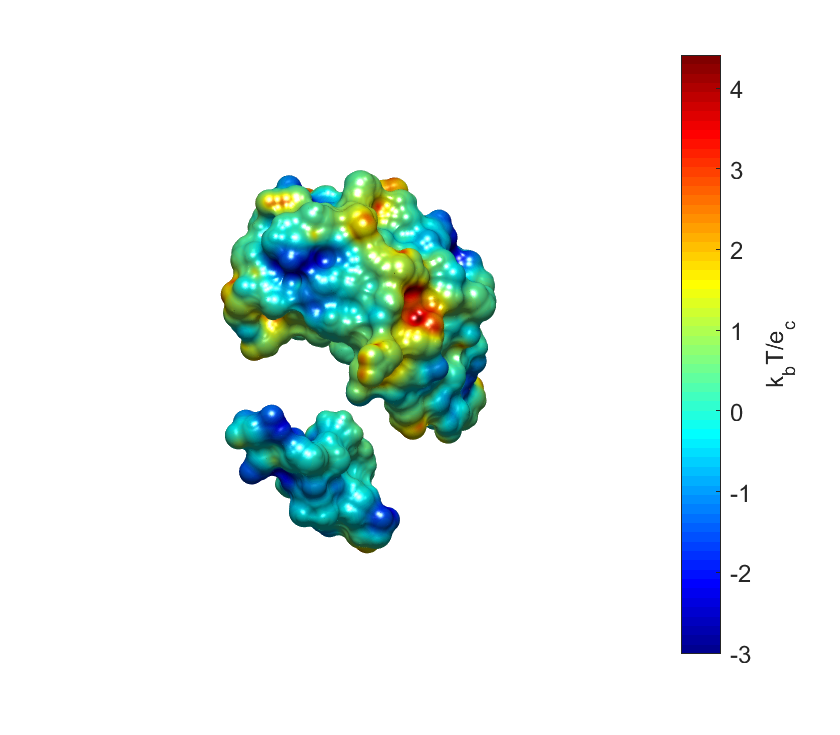}
    \caption{Electrostatic potential over the ``tight'' solvent-molecule interface. The electrostatic potential and its normal derivative $\psi$ and $\psi_n$ are approximated numerically using K-reg and CTR2, with tubular neighborhood of width $\varepsilon=2h$ and Cartesian grid of size $513^3$. The system is inverted in 13 iterations for K-reg and 15 for CTR2. Then we compute the polarization energy $\mathcal{G}_{pol}$ and obtain the following values: -1.4534e+3 for K-reg and  -1.4478e+3 for CTR2. The electrostatic potential plotted is computed using CTR2.}
    \label{fig:tight-potential}
\end{figure}

In a different test, we solve the system for the tight smooth surface by using as centers the centers of the atoms, and compute the polarization energy $\mathcal{G}_{pol}$ \eqref{eq:pol-energy}. We use both K-reg and CTR2, and see that the system is well-conditioned by comparing the number of GMRES iterations needed. The potential obtained is plotted in Figure~\ref{fig:tight-potential}. \\

\begin{table}[htb!]
    \caption{{Smooth ``loose'' solvent-molecule interface results.} 
    We solve the system on the surface presented on the left in Figure~\ref{fig:smooth-molecule-1}, using K-reg and CTR2. 
    Then we evaluate the function $\psi_{rxn}$ \eqref{eq:psi-rxn} in the origin to test the accuracy of the methods. 
    Concurrently we test the area approximation with the two different Jacobian approximations used in the two different methods.
    }\medskip
    \begin{center}
        \begin{tabular}{c|c|c|c|c|c}
            grid size & width $\varepsilon$ & d.o.f. & $\mathcal{A}$, $J_\Gamma \approx 1$ & $\psi_{rxn}(\zero)$, K-reg & GMRES \\ \hline
            $513^3$ & $2h$ & 1468195 & 5661.51 & -1.6677e+1 & 11 \\
            $513^3$ & $5h$ & 3671757 & 5662.14 & -1.6671e+1 & 11 \\[0.2cm]
            grid size & width $\varepsilon$ & d.o.f. & $\mathcal{A}$, $J_\Gamma \approx J_h$ & $\psi_{rxn}(\zero)$, CTR2 & GMRES \\ \hline
            $513^3$ & $2h$ & 1468195 & 5661.38 & -1.6684e+1 & 12  \\
            $513^3$ & $5h$ & 3671757 & 5662.12 & -1.6686e+1 & 11  
        \end{tabular}
    \end{center}
    \label{tab:smooth-surface-VISM-1}
\end{table}

\begin{table}[htb!]
    \caption{{Smooth ``tight'' solvent-molecule interface results, analogous to the ones presented in Table \ref{tab:smooth-surface-VISM-1}. The surface can be seen plotted on the right in Figure~\ref{fig:smooth-molecule-1}}.
    }\medskip
    \begin{center}
        \begin{tabular}{c|c|c|c|c|c}
            grid size & width $\varepsilon$ & d.o.f. & $\mathcal{A}$, $J_\Gamma \approx 1$ & $\psi_{rxn}(\zero)$, K-reg & GMRES \\ \hline
            $513^3$ & $2h$ & 1521649 & 5866.30 & 2.2840e+2 & 12 \\
            $513^3$ & $5h$ & 3805197 & 5867.53 & 2.2915e+2 & 12 \\[0.2cm] 
            grid size & width $\varepsilon$ & d.o.f. & $\mathcal{A}$, $J_\Gamma \approx J_h$ & $\psi_{rxn}(\zero)$, CTR2 & GMRES \\ \hline
            $513^3$ & $2h$ & 1521649 & 5866.10 & 2.3035e+2 & 13 \\ 
            $513^3$ & $5h$ & 3805197 & 5866.87 & 2.3104e+2 & 12 
        \end{tabular}
    \end{center}
    \label{tab:smooth-surface-VISM-2}
\end{table}

\subsubsection{Solvent excluded surfaces}
The next tests {involve the hybrid method, HYB-IBIM, for solvent excluded surfaces (SES), which are globally $C^1$ and piecewise $C^2$.} 

\begin{figure}[!htb]
    \begin{center}
         \includegraphics[width=0.25\textwidth]{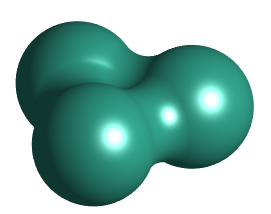}\\[0.2cm]
        \includegraphics[width=0.35\textwidth]{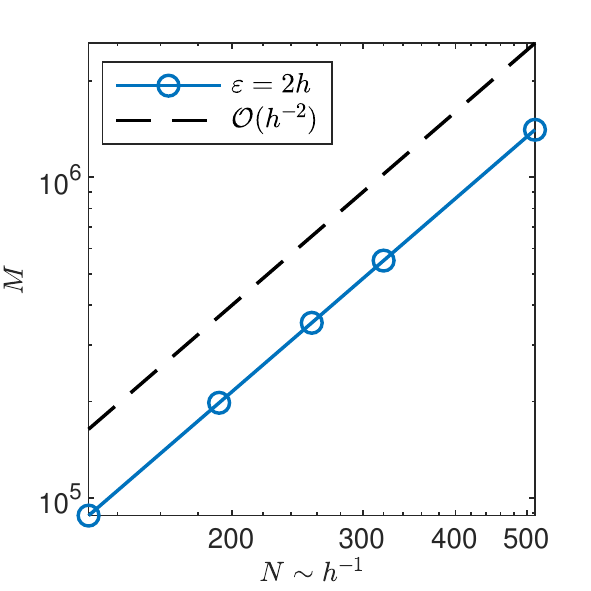}
        \includegraphics[width=0.35\textwidth]{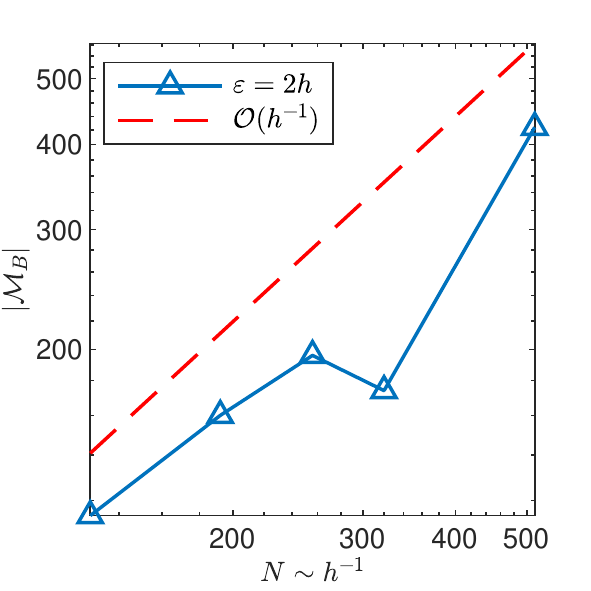}
    \end{center}
    \caption{Top: Illustration of the SES for a three-sphere model. Bottom left: $M$ number of nodes interior to the tubular neighborhood, i.e. degrees of freedom, as function of $N$ discretization along a single direction. For $\varepsilon$ linear in $h$, the number of d.o.f. scales as $N^2\sim h^{-2}$. Bottom right: number $|\mathcal{M}_B|$ of ``bad'' nodes compared with the total number of nodes scales linearly in $h^{-1}$.}
    \label{fig:three-spheres-bad-nodes}
\end{figure}

First we consider a three-sphere model shown in Figure~\ref{fig:three-spheres-bad-nodes}.
The SES is generated using the method in~\cite{zhong2018implicit}. 
{In this implementation, ``bad points'' are classified if the curvatures at that point are away from the true ones.} 
We numerically verify that the number of bad points (ref. Sec~\ref{sec: hyd cor quad}) $|\mathcal{M}_B|$ scales linearly in $h^{-1}$ if $\varepsilon\sim h$. 
as presented in~Figure~\ref{fig:three-spheres-bad-nodes}.

Next, with Table~\ref{tab:three-spheres-tab}, we present a numerical convergence study  using the three available methods: K-reg, CTR2, and HYB-IBIM. 
We compute the difference of consecutive values of $\psi_{rnx}(\zero)$ to estimate the accuracy, and see that the CTR2 and HYB methods appear to have higher accuracy than K-reg.
{In this study, it seems that the difference between CTR2 and HYB is minimal. }

\begin{table}[h!]
    \caption{Three spheres test, with width $\varepsilon=2h$. We solve the system on the surface presented in Figure~\ref{fig:three-spheres-bad-nodes}, using K-reg, CTR2, and HYB-IBIM. 
    Then we evaluate the function $\psi_{rxn}$ \eqref{eq:psi-rxn} in the origin to test the accuracy of the methods. 
    Concurrently we test the area approximation with the three different Jacobian approximations used in the three different methods. 
    The ``diff'' columns present the different of consecutive values in the column to their left. 
    The ``GMRES'' column presents the number of GMRES iterations needed to solve the system.}\medskip
    \begin{center}
        \begin{tabular}{c|c|c|c|c|c}
            grid size & $\mathcal{A}$, $J_\Gamma \approx 1$ & diff & $\psi_{rxn}(\zero)$, K-reg & diff & GMRES \\ \hline
            $128^3$ & 21.477653 & & -2.39016e+02  &  & 6 \\
            $192^3$ & 21.474494 & 3.159e-3 & -2.39906e+02 & 8.89e-1 & 7 \\
            $256^3$ & 21.473238  & 1.256e-3 &  -2.40361e+02 & 4.54e-1 & 7 \\
            $320^3$  & 21.472567 & 6.709e-4 & -2.40647e+02  & 2.86e-1 & 7 \\
            $512^3$ & 21.472018  & 5.490e-4 &  -2.41084e+02 & 4.36e-1 & 7 \\[0.2cm] 
            grid size & $\mathcal{A}$, $J_\Gamma \approx J_h$ & diff & $\psi_{rxn}(\zero)$, CTR2 & diff & GMRES \\ \hline
            $128^3$  & 21.469761 & &  -2.41738e+02  & & 10 \\
            $192^3$ & 21.471444 & 1.683e-3 &  -2.41849e+02 & 1.11e-1 & 10 \\
            $256^3$ & 21.471376  & 6.800e-5 &  -2.41844e+02 & 4.30e-3 & 10 \\
            $320^3$ & 21.471352 & 2.399e-5 &  -2.41863e+02  & 1.83e-1  & 10 \\
            $512^3$ & 21.471791  & 4.390e-4 &  -2.41881e+02  & 1.83e-2  & 10  \\[0.2cm]
            grid size & $\mathcal{A}$, $J_\Gamma \approx J_h^{(H)}$ & diff & $\psi_{rxn}(\zero)$, HYB-IBIM & diff & GMRES \\ \hline
            $128^3$ & 21.469783 & &  -2.41739e+02  &   & 10 \\
            $192^3$ & 21.471498 & 1.715e-3 &  -2.41847e+02  & 1.07e-1  & 10 \\
            $256^3$ & 21.471360  & 1.379e-4 &  -2.41842e+02 & 4.61e-3  & 10 \\
            $320^3$ & 21.471342 & 1.800e-5 &  -2.41863e+02  & 2.02e-2 & 10 \\
            $512^3$ &  21.471784  & 4.419e-4 &  -2.41880e+02 & 1.75e-2  & 10 
        \end{tabular}
    \end{center}
    \label{tab:three-spheres-tab}
\end{table}

\begin{figure}[h!]
    \centering
    \includegraphics[trim = 3cm 0cm 0cm 1.8cm, clip, scale=0.45]{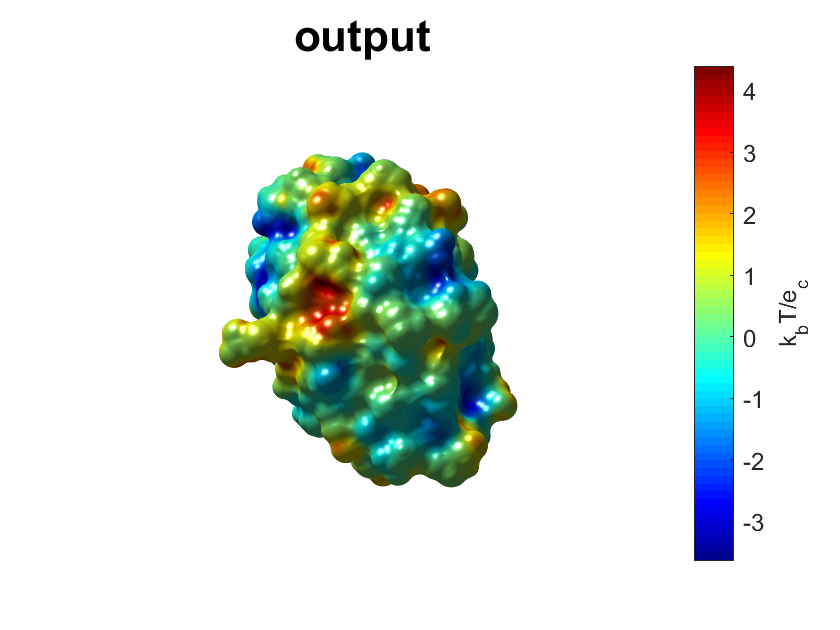}
    \caption{Electrostatic potential over the SES of the protein 1YCR, computed using HYD-IBIM.}
    \label{fig:1ycr-potential}
\end{figure}

Finally, we apply the hybrid method also on a solvent-molecule interface, we take the protein 1YCR and generate the SES using the method from \cite{zhong2018implicit}. We run a convergence test on it, seen in Table~\ref{tab:protein-SES-1YCR}. The electrostatic potential computed with the HYB-IBIM is also plotted for $513^3$ in Figure~\ref{fig:1ycr-potential}. 

\begin{table}[h!]
    \caption{Protein 1YCR. We solve the system on the SES of protein 1YCR using K-reg, CTR2, and HYB-IBIM, with width $\varepsilon=2h$. 
    Then we compute the polarization energy $\mathcal{G}_{pol}$ \eqref{eq:pol-energy} to test the accuracy of the methods.}\medskip
    \begin{center}
        \begin{tabular}{c|c|c|c|c|c}
            grid size & $\mathcal{A}$, $J_\Gamma \approx 1$ & diff & $\mathcal{G}_{pol}$, K-reg & diff & GMRES \\ \hline
            $128^3$ & 4388.871 & & -1.1905865e+03  &  & 11 \\
            $256^3$ & 4461.273  & 7.24e+1 & -1.1951553e+03 & 4.56e+0 & 12 \\
            $512^3$ & 4500.038  & 3.87e+1 & -1.2025009e+03 & 7.34e+0 & 11\\[0.2cm] 
            grid size & $\mathcal{A}$, $J_\Gamma \approx J_h$ & diff & $\mathcal{G}_{pol}$, CTR2 & diff & GMRES \\ \hline
            $128^3$ & 4366.619 &  & -1.1584408e+03 & & 15 \\
            $256^3$ & 4456.909 & 9.02e+1 & -1.1824920e+03 & 2.40e+1 & 15 \\
            $512^3$ & 4499.230 & 4.23e+1 & -1.1973823e+03 & 1.48e+1 & 15  \\[0.2cm]
            grid size & $\mathcal{A}$, $J_\Gamma \approx J_h^{(H)}$ & diff & $\mathcal{G}_{pol}$, HYB-IBIM & diff & GMRES \\ \hline
            $128^3$ & 4372.020 &  & -1.1602900e+03 &   & 14 \\
            $256^3$ & 4457.470 & 8.54e+1 & -1.1826365e+03 & 2.23e+1 & 14 \\
            $512^3$ & 4499.400 & 4.19e+1 & -1.1974421e+03 & 1.48e+1 & 15 
        \end{tabular}
    \end{center}
    \label{tab:protein-SES-1YCR}
\end{table}

\section{Conclusions}

In this paper, we present an accurate corrected trapezoidal-rule based the numerical solution of the linearized Poisson-Bolztmann equation. 
Using the boundary integral formulation developed in Juffer \emph{et al.} \cite{juffer1991electric} we solve the boundary integral equations using a non-parametric approach. 
In theory \cite{izzo2022convergence} it is possible to derive very high order accurate quadratures for corresponding boundary integral equations. 
However, high order accuracy requires  higher order accurate approximations to the surface's geometry.
We identify the kernel $K_{21}$ as the bottleneck for the order of accuracy attainable. 
Given a second order approximation of the surface locally around the target points, the {proposed corrected trapezoidal rule can} reach third order of accuracy for integrating all the kernels {except} $K_{21}$.  
To obtain a third order approximation for integrating $K_{21}$, a fourth order approximation of the surface geometry is needed; see Section~\ref{subsec:expansions}. 

We first tested the method on a single sphere/ion and checked the order of accuracy of the method. Moreover we assessed that the error constant when using a width of $\varepsilon$ linear in $h$ is much smaller than the regularization method K-reg.

Then we studied the application of the Poisson-Boltzmann equation to compute the electrostatic potential and polarization energy of macromolecules immersed in a solvent. The smoothness of the solvent-molecule interface depends on how it is defined. We presented polarization energies for the smooth interfaces, and then tested a hybrid corrected/regularized method for piecewise $C^2$ surfaces.

\section*{Acknowledgments}

The authors would like to thank Zirui Zhang and Prof. Li-Tien Cheng for providing the smooth solvent-molecule interface for the p53-MDM2 protein via its signed distance function. 
The authors would also like to thank the Texas Advanced Computing Center (TACC) and the Swedish National Infrastructure for Computing at PDC Center for High Performance Computing for the computing resources. 
Tsai’s research is supported partially by NSF grant DMS-2110895.

\printbibliography
\end{document}